\documentclass[11pt]{article}
\usepackage[cm]{fullpage}
\usepackage{amscd}
\usepackage{amsmath}
\usepackage{amsfonts}
\newcommand{\al}{\alpha}
\newcommand{\be}{\beta}
\newcommand{\ga}{\gamma}
\newcommand{\de}{\delta}

\newcommand{\dd}[1]{\partial_{#1}}

\newcommand{\ov}[1]{\overline{#1}}
\newcommand{\un}[1]{\underline{#1}}
\newcommand{\J}[2]{J_{#1}^{\ #2}}
\newcommand{\PP}[4]{\mathcal{P}^{#1 #2}_{#3 #4}}
\newcommand{\QQ}[4]{\mathcal{Q}^{#1 #2}_{#3 #4}}

\newcommand{\tr}[2]{\textrm{tr}_{#1}#2}
\newcommand{\G}[2]{G_{#1 \overline{#2}}}
\newcommand{\GI}[2]{G^{#1 \overline{#2}}}
\newcommand{\Gp}[2]{G'_{#1 \overline{#2}}}
\newcommand{\GpI}[2]{G'^{#1 \overline{#2}}}
\newcommand{\GT}[2]{\tilde{G}_{#1 \overline{#2}}}
\newcommand{\GTI}[2]{\tilde{G}^{#1 \overline{#2}}}
\newcommand{\GTp}[2]{\tilde{G}'_{#1 \overline{#2}}}

\newcommand{\ddbar}{\sqrt{-1} \partial \ov{\partial}}
\newcommand{\ot}{\tilde{\omega}}
\newcommand{\Rm}{\textrm{Rm}}

\begin{document}
\newcounter{remark}
\newcounter{theor}
\setcounter{remark}{0}
\setcounter{theor}{1}
\newtheorem{claim}{Claim}
\newtheorem{theorem}{Theorem}[section]
\newtheorem{proposition}{Proposition}[section]
\newtheorem{lemma}{Lemma}[section]
\newtheorem{defn}{Definition}[theor]
\newtheorem{corollary}{Corollary}[section]
\newenvironment{proof}[1][Proof]{\begin{trivlist}
\item[\hskip \labelsep {\bfseries #1}]}{\end{trivlist}}
\newenvironment{remark}[1][Remark]{\addtocounter{remark}{1} \begin{trivlist}
\item[\hskip
\labelsep {\bfseries #1  \thesection.\theremark}]}{\end{trivlist}}

\centerline{\bf \Large THE CALABI-YAU EQUATION}
\centerline{\bf \Large ON ALMOST-K\"AHLER FOUR-MANIFOLDS}
\bigskip
\bigskip

\centerline{\bf Ben Weinkove\footnote{The author is supported in part by National Science Foundation grant DMS-05-04285.  This work was carried out while the author was visiting Imperial College, London on a Royal Society Research Assistantship.}}
\centerline{Harvard University}
\centerline{Department of Mathematics}
\centerline{Cambridge, MA 02138} 
\bigskip
\bigskip
\bigskip

\noindent
{\bf Abstract.} Let $(M, \omega)$ be a compact symplectic 4-manifold with
a compatible almost complex structure $J$.  The problem
of finding a $J$-compatible symplectic form
 with prescribed volume form is an almost-K\"ahler analogue of Yau's
 theorem and is connected to a programme in symplectic topology proposed
 by Donaldson.   We call the corresponding  equation for the symplectic form the \emph{Calabi-Yau equation}.
  Solutions are unique in their cohomology class.    It is shown in this
  paper that a solution to this equation
  exists if the Nijenhuis tensor is small in a certain sense.  Without this assumption, it is shown
 that the problem of existence can be
 reduced to obtaining  a $C^0$ bound on  a scalar
 potential function.

\bigskip
\bigskip

\setlength\arraycolsep{2pt}
\addtocounter{section}{1}
\noindent
{\bf 1. Introduction}
\bigskip

In 1954 Calabi \cite{Ca}  conjectured that any representative of the first
Chern class of a compact K\"ahler manifold $(M, \omega)$ can be written as the Ricci
curvature of a K\"ahler metric $\omega'$ cohomologous to $\omega$.  He showed that any such
metrics are unique.  Yau \cite{Ya} famously
solved Calabi's conjecture
around twenty years later.  This result, and the immediate corollary that any K\"ahler
manifold with $c_1(M)=0$ admits a Ricci-flat metric, have had many applications in both mathematics and theoretical physics.  

Yau's theorem is equivalent to finding a K\"ahler metric in a given K\"ahler class with prescribed volume form.  By the $\partial \overline{\partial}$-Lemma this amounts to
solving the complex Monge-Amp\`ere equation
\begin{equation} \label{eqnCY1}
(\omega + \ddbar \phi)^n = e^F \omega^n,
\end{equation}
for smooth real $\phi$ with $\omega + \ddbar \phi>0$, where  $n=\textrm{dim}_{\mathbb{C}} M$ and
  $F$
is any smooth function with $\int_M e^F \omega^n = \int_M \omega^n$.
Yau solved this equation  by considering the family of equations obtained
by replacing $F$ by $tF+c_t$
for some constant $c_t$, for $t \in [0,1]$ and using the continuity method.  This requires an
openness argument using the implicit function theorem, and, more importantly,
a  closedness argument which requires his celebrated \emph{a priori} estimates.
Yau also generalized Calabi's conjecture: first in the case when the right hand
side of (\ref{eqnCY1}) may have poles or zeros \cite{Ya}; and second, with Tian, in the context of complete non-compact K\"ahler Ricci-flat metrics 
\cite{TiYa1, TiYa2}.  For other results along these lines, see  \cite{Ko}, \cite{BaKo1}, \cite{BaKo2}, \cite{Jo}, for example.

The aim of this paper is to attempt to generalize  Yau's theorem in a very different
direction.  We consider the 
 case when  
the almost complex structure is not integrable.  This problem was suggested to the
author by Donaldson  and is motivated by a wider programme of his on the  symplectic
topology of 4-manifolds \cite{Do}.
 Let $(M,\omega)$ be a symplectic
four-manifold.  Then there exists an almost complex structure $J$ which is
compatible with $\omega$.  This defines a metric $g$ by
$$g( \cdot, \cdot) = \omega ( \cdot, J \cdot )>0.$$
If $J$ is integrable then it is K\"ahler.  In general, the data $(M, \omega, J)$ is called an \emph{almost-K\"ahler
manifold} and we will call $\omega$ an \emph{almost-K\"ahler form}.  The volume form version of Yau's theorem still makes sense.
Given an almost-K\"ahler 4-manifold $(M, \omega, J)$  we ask whether there exists an almost-K\"ahler form $\omega'$  solving the equation
\begin{equation} \label{eqnCY2}
{\omega'}^2 = e^F \omega^2,
\end{equation} 
for any function $F$ satisfying
\begin{equation} \label{eqnconditionF}
\int_M e^F \omega^2 = \int_M \omega^2,  
\end{equation}
and we  also ask whether $\omega'$ can be taken to be cohomologous to $\omega$.
We call  (\ref{eqnCY2})
the \emph{Calabi-Yau equation}.  Any solution
to (\ref{eqnCY2}) is unique in its cohomology class  -
this fact was pointed out to the author by Donaldson. A proof is given in
section 2.

Following Yau, we use the continuity method to try to obtain the existence of a solution.
First, we consider  the question of \emph{a priori} estimates for solutions to (\ref{eqnCY2}).  For simplicity,  assume that $\omega'$ is cohomologous to $\omega$.  We show that all the estimates can be reduced
to a uniform bound of a scalar potential
function $\phi_1$ defined, up to a constant, by
\begin{equation} \label{eqnphidefn}
\frac{\omega \wedge \omega'}{\omega'^2} = 1 - \frac{1}{4} \Delta' \phi_1,
\end{equation}
where $\Delta'$ is the Laplacian associated to  $\omega'$. The function
$\phi_1$ belongs to a   
 a  1-parameter family
of `almost-K\"ahler potentials' $\{ \phi_s\}_{s \in [0,1]}$, defined in section
2, which all coincide
in the K\"ahler case with the usual
K\"ahler potential.

\bigskip
\noindent
{\bf Theorem 1} \,  \emph{Let $(M,\omega,J)$ be a compact almost-K\"ahler 4-manifold.
Suppose that $\omega'$ is another almost-K\"ahler form, cohomologous to $\omega$, and satisfying (\ref{eqnCY2}).   Then there exist positive constants $K_i$ depending only on $(M,\omega,J)$, $F$ and $\emph{osc}_M \phi_1$ such that $\omega' \ge K_0^{-1} \omega$ and
$$ \| \omega' \|_{C^i (g)} \le K_i \quad \textrm{for } i=0,1,2,\ldots, $$
where  $\emph{osc}_M \phi_1 = \sup_M \phi_1 - \inf_M
\phi_1$. }
\bigskip
 
An analogous result holds even if 
 $\omega'$  and $\omega$ are not necessarily cohomologous (see section 7).
We turn now to the question of openness in the continuity method.  Denote
by $\mathcal{H}^+_{\omega}$ the space of self-dual harmonic 2-forms with
respect to $\omega$ and by $H^+_{\omega}$ the corresponding 
subspace in $H^2(M; \mathbb{R})$.  $H^+_{\omega}$ is a maximal positive subspace for the intersection form on $H^2(M; \mathbb{R})$ and its dimension is  $b^+(M)$.  
Notice that $\omega$ is harmonic and self-dual and so $b^+(M)\ge 1$.  If $b^+(M)=1$
 then we can show the openness part of the continuity method, remaining in
 the same cohomology class.  In the case where $b^+(M)>1$  the openness argument still works
if we allow the  class to vary within $H^+_{\omega}$.

Under  the  assumption that the Nijenhuis tensor $N(J)$ is small in the $L^1$ norm, the required uniform bound on $\phi_1$ can be obtained.
 So in this case, we can solve   equation (\ref{eqnCY2}).  

\bigskip
\pagebreak[3]
\noindent
{\bf Theorem 2} \, \emph{Let $(M,\omega, J)$ be a compact almost K\"ahler
 4-manifold. 
\begin{enumerate}
\item[(i)]  Suppose $b^+(M)=1$.  Then for $F \in C^{\infty}(M)$ satisfying (\ref{eqnconditionF})
there exists an almost-K\"ahler form
$\omega'$ cohomologous to $\omega$ solving (\ref{eqnCY2}) if
\begin{equation} \label{eqnNJ}
\| N(J) \|_{L^1(g)}< \epsilon, 
\end{equation}
for $\epsilon>0$ depending only on $g$ and $\| F\|_{C^2(g)}$.
\item[(ii)] If $b^+(M)>1$ then the same holds except that the solution $\omega'$ may lie in a different cohomology class in $H^+_{\omega}$.
\end{enumerate}}

With a little work, an explicit $\epsilon$ could be written down.  
However, it is hoped that
the condition (\ref{eqnNJ}) could be removed entirely (cf. \cite{Do}).
In addition, it would be interesting to improve on Theorem 2 even further in light
of a possible application in symplectic topology described to the author by Donaldson.   Given an almost complex
structure $J_1$ on a symplectic 4-manifold,  a natural question is:  does there exist 
a symplectic form compatible
with $J_1$?   In general the answer is negative, as can be seen from the
well-known Kodaira-Thurston example \cite{Th}, \cite{FeGoGr}.  However, it is  sensible to ask this question  
under the (obviously necessary) assumption that there exists some symplectic
form $\Omega$ taming $J_1$. In this case, there exists an almost complex structure $J_0$  compatible
with $\Omega$ and, by a well-known result of Gromov \cite{Gr}, a smooth path of almost complex structures $\{ J_t \}_{t\in[0,1]}$  all taming $\Omega$.  Set $\omega_0 = \Omega$ and consider
the equation
$$\omega_t^2 = \Omega^2,$$
for $\omega_t$ compatible with $J_t$.   Finding a solution
for $t=1$ would solve the problem.  To prove this using a continuity method  one would require  estimates for $\omega_t$ depending only on $\Omega$
and $J_t$.

These methods appear to make sense only in four dimensions, since
 the system of equations is overdetermined in higher dimensions.  Nevertheless, it should
be noted that
many of the estimates here carry over easily to any dimension.

The outline of the paper is as follows:  in section 2, some preliminaries are given, almost-K\"ahler potentials are defined and uniqueness for the Calabi-Yau equation is proved;  in section 3  an estimate
on the metric $g'$ in terms of the potential is given; a H\"older estimate
on the metric is proved in section 4; the higher order estimates and the
proof of Theorem 1 are given in section 5;  finally, in sections 6 and 7, Theorem
2 is proved in the cases $b^+(M)=1$ and $b^+(M)>1$ respectively.

\bigskip
\noindent
{\bf Remark 1.1} \ Delano\"e \cite{De} considered, following a suggestion of Gromov, a different problem concerning the equation (\ref{eqnCY2}).
He looked for solutions
of $\omega'^n = e^F \omega^n$, on an almost-K\"ahler manifold $(M, \omega)$ of dimension
$2n$,  of the form
$\omega' = \omega + d(Jd\phi)$ for a smooth real function $\phi$ so that $\omega'$
tames $J$ but is not necessarily compatible with $J$ (here, $J$ acts on 1-forms in the usual way).  He showed that when $n=2$, if  there exists such a solution for every $F$, then $J$ is in fact integrable. We
do not expect the
solutions we obtain in Theorem 2 to be, in general, of the form $\omega' = \omega + d(Jd\phi)$
for any $\phi$.

\pagebreak[3]
\bigskip
\bigskip
\addtocounter{section}{1}
\setcounter{equation}{0}
\noindent
{\bf 2. Almost-K\"ahler geometry and the Calabi-Yau equation}

\bigskip
\noindent
{\it Notation and preliminaries}
\bigskip

We will often work in local coordinates, making use of the Einstein summation convention.  The almost complex
structure $J
= \J{i}{j} dx^i \otimes \frac{\partial}{\partial x^j}$ satisfies, by definition,
the condition
$$\J{i}{k} \J{k}{j} = - \delta_{i}^{j}.$$
We will lower indices in the usual way using the metric $g$ so that
$$J_{ij} = \J{i}{k} g_{kj} = \omega_{ij}.$$
The condition $d\omega=0$ can be written as
\begin{equation} \label{eqndomega}
\dd{i} J_{jk} + \dd{j} J_{ki} + \dd{k} J_{ij}=0.
\end{equation}
It follows that the equation
\begin{equation} \label{eqnharmonic}
\nabla_i \J{j}{i}=0
\end{equation}
holds on an almost-K\"ahler manifold, where $\nabla$ is the Levi-Civita connection
associated to the metric $g$.  This implies that $\omega$ is harmonic with
respect to the metric $g$.

Define two tensors $\PP{}{}{}{}$ and $\QQ{}{}{}{}$ by 
\begin{eqnarray*}
\PP{i}{j}{k}{l} & = & \frac{1}{2} (\delta_{k}^i \delta_l^j - \J{k}{i} \J{l}{j})
\\
\QQ{i}{j}{k}{l} & = & \frac{1}{2} (\delta_{k}^i \delta_l^j + \J{k}{i} \J{l}{j}).
\end{eqnarray*}
Then the compatibility of $g$ and $\omega$ with $J$ implies that
$\PP{i}{j}{k}{l} g_{ij} = 0$ and $\PP{i}{j}{k}{l} J_{ij} = 0$.
Considering  $\PP{}{}{}{}$ and $\QQ{}{}{}{}$ as operators on two-tensors,
we have
$$\PP{}{}{}{} + \QQ{}{}{}{} = \textrm{Id}.$$
Moreover, at each point,  $\PP{}{}{}{}$ and $\QQ{}{}{}{}$ are self-adjoint with respect
to $g$ and define projections onto the spaces $\ker
\QQ{}{}{}{}$ and $\ker \PP{}{}{}{}$ respectively.

The obstruction to the almost complex structure $J$ being integrable is  the Nijenhuis tensor $N: TM \times TM \rightarrow TM$, which is given by
$$N(X,Y) = [X,Y] + J [JX, Y] + J[X, JY] - [JX,JY] .$$
In local coordinates, this can be written as
$$N_{jk}^i =  \J{k}{l} \dd{l} \J{j}{i} +  \J{l}{i} \dd{j} \J{k}{l} - \J{j}{l} \dd{l} \J{k}{i}  
- \J{l}{i} \dd{k} \J{j}{l}.$$
On an almost-K\"ahler manifold, the Nijenhuis tensor can be written in the
simpler form
$$N_{jk}^i = 2  (\nabla^i \J{j}{l} ) J_{kl}.$$
By the Newlander-Nirenberg theorem, the almost complex structure $J$ is integrable if and only if $N$ vanishes
identically and if and only if $\nabla J$ vanishes identically.

For later use, we also make the following simple observation.  Let $\Rm$ denote
the Riemmanian curvature tensor of an almost-K\"ahler metric $g$.  Then
\begin{equation} \nonumber
\sup_M |\nabla J |^2 \le C \| \Rm \|_{C^0(g)},
 \end{equation}
for $C$ a constant depending only on dimension. Indeed, using the usual commutation formulae
for covariant derivatives along with (\ref{eqndomega}) and (\ref{eqnharmonic}),
$$0 = \Delta | J|^2 = 2 | \nabla J |^2 + \Rm * J * J,$$
where $*$ denotes some bilinear operation involving tensor products and the
metric $g$.  Similarly, by calculating $\Delta |\nabla J|^2$ we see that
$$ | \nabla \nabla J|^2 = \nabla \nabla \Rm * J * J + \Rm * \nabla \nabla J *
J + \nabla \Rm * J * \nabla J + \Rm * \nabla J * \nabla J,$$
and it follows that $\sup_M|\nabla \nabla J|^2$ can be bounded by a constant depending
only on dimension and $\| \Rm \|_{C^2(g)}$.

\bigskip
\pagebreak[3]
\noindent
{\it Almost-K\"ahler potentials}
\bigskip

Now restrict to four dimensions.
Let $\omega$ and $\omega'$ be two almost-K\"ahler forms with $[\omega] = [\omega']$.  For $s \in [0,1]$, let  $\Omega_s = (1-s) \omega + s \omega'$ and define the \emph{almost-K\"ahler potentials} $\phi_s$ by
$$(1-2s) \omega \wedge \omega' + s \omega'^2 - (1-s)\omega^2 = - \frac{1}{2} \Omega_s \wedge d ( J d \phi_s).$$
Since $- \Omega \wedge d (J d\phi)
= \frac{1}{2} \Delta_{\Omega} \phi \, \Omega^2$ for any almost-K\"ahler form $\Omega$ and function
$\phi$, the  existence of the $\phi_s$ follows from elementary Hodge theory.
 The $\phi_s$ are uniquely determined up to the addition of a constant.
 In the K\"ahler case, they all coincide with the usual K\"ahler potential
$\phi$ given by $\omega' = \omega + \sqrt{-1} \partial \ov{\partial} \phi$. We are interested in three particular almost-K\"ahler potentials, corresponding to $s=0, \frac{1}{2}, 1$.  They satisfy:
\begin{eqnarray}  
\frac{1}{4} \Delta \phi_0 & = & \frac{\omega \wedge \omega'}{\omega^2} - 1 \\
\frac{1}{4} \Delta' \phi_1 & = & 1 -  \frac{\omega \wedge \omega'}{\omega'^2}\\
\label{eqnOmega}
\frac{1}{2} (\Delta_{\frac{1}{2}} \phi_{\frac{1}{2}})  & = & \frac{\omega'^2 - \omega^2}{\Omega_{\frac{1}{2}}^2},  
\end{eqnarray}
where $\Delta'$ and $\Delta_{\frac{1}{2}}$ are the Laplacians associated to $\omega'$ and $\Omega_{\frac{1}{2}}$.  

In addition, for each $s$, define a one form $a_s$ by the equations
$$\omega' = \omega - \frac{1}{2}d( J d\phi_s) + da_s,$$
and $d^*_s a_s=0$, where $d^*_s$ is the formal adjoint of $d$ associated to the metric $\Omega_s$.  Note that $a_s$ is defined only up to the addition of a harmonic 1-form.  A short calculation shows that $a_s$ satisfies the elliptic system
\begin{equation} \label{eqnellipticsystem1}
\left. \begin{array}{rcl} 
da_s \wedge \Omega_s & = & 0 \\
\PP{}{}{}{} da_s & = & \frac{1}{4} (\dd{i} \J{j}{k} - \dd{j} \J{i}{k}) (\dd{k} \phi_s) \, dx^i \wedge dx^j  \\
d^*_s a_s & = & 0. \end{array} \right\}
\end{equation}
It will be convenient to give a different formulation of (\ref{eqnellipticsystem1}).  Let $*_s$ be the Hodge-star operator associated to $\Omega_s$.  Then the projection $\frac{1}{2} (1 + *_s) : \Lambda^2 \rightarrow \Lambda_s^+$ onto the self-dual two forms can be written
$$\frac{1}{2} (1+ *_s) (\chi) = \left( \frac{\Omega_s \wedge \chi}{\Omega_s^2} \right) \Omega_s + \PP{}{}{}{} \chi, \quad \textrm{for } \chi \in \Lambda^2.$$
Hence, the system (\ref{eqnellipticsystem1}) can be rewritten as
\begin{equation} \label{eqnellipticsystem1.5}
\left. \begin{array}{rcl} 
d^+_s a_s & = & \frac{1}{4} (\dd{i} \J{j}{k} - \dd{j} \J{i}{k}) (\dd{k} \phi_s) \, dx^i \wedge dx^j \\
d^*_s a_s & = & 0,
\end{array} \right\}
\end{equation}
for $d^+_s = \frac{1}{2} (1+*_s) d : \Lambda^1 \rightarrow \Lambda_s^+$.
 Observe that, since $d_s^+b=0$ implies $db=0$ for any 1-form $b$ (see for example \cite{DoKr}, Prop. 1.1.19), the kernel of the operator $(d_s^+, d_s^*)$ consists of the harmonic 1-forms.

\bigskip
\pagebreak[3]
\noindent
\emph{The Calabi-Yau equation}
\bigskip

For a manifold of dimension $n=2m$,  linearizing the Calabi-Yau equation gives
$a \mapsto {\omega^{m-1} \wedge da}$ for a one-form $a$.  Combining with the linear operator $a \mapsto \PP{}{}{}{} da$ and imposing  $d^*a=0$ gives a system which is elliptic if $n=4$ and overdetermined if $n >  4$.

Restricting to four-manifolds, 
solutions to (\ref{eqnCY2}) are unique in their cohomology class.  To see this, let $\omega_1$ and
$\omega_2= \omega_1+db$ be cohomologous almost-K\"ahler forms with $\omega_1^2=\omega_2^2$.
Let $\Omega = \omega_1 + \omega_2$. Then $\Omega \wedge db = 0$ and $\PP{}{}{}{}db=0$
from which it follows that $d^+_{\Omega}b=0$.  Then $db=0$ and so $\omega_1
= \omega_2$.

Finally we mention that the Calabi-Yau equation (\ref{eqnCY2}) can also be written as
\begin{equation} \label{eqnCY3}
\tr{g}{g'} = e^F \tr{g'}{g},
\end{equation}
where we are writing $g'$ for the metric associated to $\omega'$ and where
$$ \tr{g}{g'} = g^{ij} g'_{ij} = J^{ij} J'_{ij}, \quad \textrm{and} \quad \tr{g'}{g}=  g'^{ij}g_{ij} = J'^{ij} J_{ij}.$$
We are using here the obvious notation for lowering and raising indices using the metric $g'$, so that $J'_{ij} = \J{i}{k} g'_{kj}$ and
$J'^{ij} = \J{k}{j}g'^{ki}$.


\setcounter{equation}{0}
\addtocounter{section}{1}
\bigskip
\bigskip
\pagebreak[3]
\noindent
{\bf 3. Uniform estimate on the metric}

\bigskip

In this section we will prove an estimate on the metric in terms of $\phi_1$,
similar to the one proved by Yau \cite{Ya} (see also \cite{Au}) in the
K\"ahler case.  The computations here are somewhat more involved because
of extra terms arising from the non-integrability of the almost complex structure.

\begin{theorem} \label{theoremgestimate}
There exist constants $C$ and $A$ depending only on  $\| \emph{Rm}(g) \|_{C^2(g)}$,
$\sup_M |F|$ and the lower
bound of $\Delta F$  such that
$$\emph{tr}_g{g'} \le C e^{A(\phi_1 - \inf_{M} \phi_1)}.$$
\end{theorem}
\begin{proof}
We will calculate $\Delta' (\log (\tr{g}{g'}) -A \phi_1)$, where $A$ is a constant to be determined, and then apply the maximum principle. For this
calculation we will work at a point  using normal coordinates for $g$.

First, the equation (\ref{eqnCY2}) can be written
\begin{equation} \label{eqnCY4}
\log \det g'  = 2F + \log \det g.
\end{equation}
Applying the Laplace operator $\Delta$ of the metric $g$, we obtain
\begin{equation} \label{eqnlapF}
2\Delta F = g^{ij} g'^{kl} \dd{i} \dd{j} g'_{kl} - g^{ij} g'^{kq} g'^{pl} ( \dd{i} g'_{pq} ) (\dd{j} g'_{kl}) - g^{ij} g^{kl} \dd{i} \dd{j} g_{kl}.
\end{equation}
Now calculate
\begin{eqnarray*}
\Delta' (\tr{g}{g'}) & = &g'^{kl} g^{ij} \dd{k} \dd{l} g'_{ij} + g'^{kl} g'_{ij} \dd{k} \dd{l} g^{ij} - g'^{kl} {\Gamma'}_{kl}^p g^{ij} \dd{p} g'_{ij},
\end{eqnarray*}
where we use ${\Gamma'}_{kl}^p$ to denote the Christoffel symbols corresponding to $\nabla'$.
We can rewrite this last term using the equation (valid at a point)
$$g'^{kl} {\Gamma'}^p_{kl} = - \J{q}{p} g'^{kl} \nabla_k \J{l}{q},$$
which follows from the equation $g'^{kl} \nabla'_k \J{l}{q} = 0$.
Thus we have
\begin{equation} \label{eqnlaptrgg'}
\Delta' (\tr{g}{g'})  = g'^{kl} g^{ij} \dd{k} \dd{l} g'_{ij} + g'^{kl} g'_{ij} \dd{k} \dd{l} g^{ij} + \J{q}{p} g'^{kl} ( \nabla_k \J{l}{q}) g^{ij} \dd{p} g'_{ij}.
\end{equation}
The first terms on the right hand side of (\ref{eqnlapF}) and (\ref{eqnlaptrgg'}) are related as follows:
\begin{eqnarray} \nonumber
g'^{kl} g^{ij} \dd{i} \dd{j} g'_{kl} & = & g'^{kl} g^{ij} \dd{k} \dd{l} g'_{ij} - \J{r}{k} g^{ij} \dd{i} \dd{j} \J{k}{r} - 2J'_{jr} g'^{kl} g^{ij} \dd{i} \dd{l} \J{k}{r} \\ \nonumber
&& \mbox{} + J'_{rj} g'^{kl} g^{ij} \dd{k} \dd{l} \J{i}{r}  - 2 \J{r}{k} g^{ij} \dd{k} \dd{j} \J{i}{r} \\ \nonumber
&& \mbox{} -  
2g'^{kl} g^{ij} ( \dd{i} \J{k}{r} + \dd{k} \J{i}{r}) \J{j}{q} \dd{r} g'_{ql} \\ \nonumber
&& \mbox{} - 2g^{ij} (\dd{i} \J{q}{r} + \dd{q} \J{i}{r} ) \dd{r} \J{j}{q} \\ \label{eqnddg}
&& \mbox{} - 2g'^{kl} g^{ij} ( \dd{l} \J{k}{r} ) \J{j}{s} \dd{i} g'_{rs}.
\end{eqnarray}
To see this, calculate
\begin{eqnarray} \nonumber
g'^{kl} g^{ij} \dd{i} \dd{j} g'_{kl} & = & - g'^{kl} g^{ij} \dd{i} \dd{j} (\J{k}{r} J'_{rl}) \\ \nonumber
& = & - \J{k}{r} g'^{kl} g^{ij} \dd{i} \dd{j} J'_{rl} - J'_{rl} g'^{kl} g^{ij} \dd{i} \dd{j} \J{k}{r} \\ \label{eqnddg2}
&& \mbox{} - 2g'^{kl} g^{ij} (\dd{i}  \J{k}{r}) (\dd{j} J'_{rl}).
\end{eqnarray}
We will now apply (\ref{eqndomega}) to the first term on the right hand side of the above equation to obtain
\begin{eqnarray*}
- \J{k}{r} g'^{kl} g^{ij} \dd{i} \dd{j} J'_{rl} & = &- \J{k}{r} g'^{kl} g^{ij} \dd{i} (\dd{r} J'_{jl} + \dd{l} J'_{rj}) \\
& = & J'^{lr} g^{ij} (\dd{i} \dd{l} J'_{jr} + \dd{i} \dd{l} J'_{jr}) \\
& = & 2 \J{k}{r} g'^{kl} g^{ij} \dd{i} \dd{l} J'_{jr} \\
& = & 2 g'^{kl} g^{ij} \dd{i} \dd{l} (\J{k}{r} J'_{jr}) - 2g'^{kl} g^{ij} (\dd{i} \J{k}{r} ) (\dd{l} J'_{jr}) \\
&& \mbox{} - 2 g'^{kl} g^{ij} (\dd{l} \J{k}{r}) (\dd{i} J'_{jr}) - 2J'_{jr} g'^{kl} g^{ij} \dd{i} \dd{l} \J{k}{r}.
\end{eqnarray*}
 \end{proof}
Then in (\ref{eqnddg2}) we have
\begin{eqnarray} \nonumber
g'^{kl} g^{ij} \dd{i} \dd{j} g'_{kl} & = & 2 g'^{kl} g^{ij} \dd{i} \dd{l} g'_{kj} - 2g'^{kl} g^{ij} (\dd{i} \J{k}{r} ) (\dd{l} J'_{jr}) \\ \nonumber
&& \mbox{} - 2 g'^{kl} g^{ij} (\dd{l} \J{k}{r}) (\dd{i} J'_{jr}) - 2J'_{jr} g'^{kl} g^{ij} \dd{i} \dd{l} \J{k}{r} \\ \label{eqnddg3} 
&& \mbox{}  - J'_{rl} g'^{kl} g^{ij} \dd{i} \dd{j} \J{k}{r} - 2g'^{kl} g^{ij} (\dd{i}  \J{k}{r}) (\dd{j} J'_{rl}). 
\end{eqnarray}
Similarly,
\begin{eqnarray} \nonumber
g'^{kl} g^{ij} \dd{k} \dd{l} g'_{ij} & = & 2 g^{ij} g'^{kl} \dd{k} \dd{j} g'_{il} - 2g^{ij} g'^{kl} (\dd{k} \J{i}{r}) (\dd{j} J'_{lr}) \\ \nonumber
&& \mbox{} - 2g^{ij} g'^{kl} (\dd{j} \J{i}{r} ) (\dd{k} J'_{lr}) - 2J'_{lr} g^{ij} g'^{kl} \dd{k} \dd{j} \J{i}{r} \\ \label{eqnddg4}
&& \mbox{} - J'_{rj} g'^{kl} g^{ij} \dd{k} \dd{l} \J{i}{r} - 2 g'^{kl} g^{ij} (\dd{k} \J{i}{r}) (\dd{l} J'_{rj}).
\end{eqnarray}
Combining (\ref{eqnddg3}) and (\ref{eqnddg4}) and again making use of (\ref{eqndomega}) gives (\ref{eqnddg}).
From (\ref{eqnlapF}), (\ref{eqnlaptrgg'}) and (\ref{eqnddg}) we obtain
\begin{eqnarray} \nonumber
\Delta' (\log \tr{g}{g'}) & = & \frac{1}{\tr{g}{g'}} \left( \Delta' (\tr{g}{g'})
- \frac{|\nabla' \tr{g}{g'}|^2}{\tr{g}{g'}} \right) \\ \nonumber
& = & \frac{1}{\tr{g}{g'}} \left( 2\Delta F + g'^{kl} g'^{ij} g^{pq} (\nabla_p g'_{ik}) ( \nabla_q g'_{jl})  \right. \\ \nonumber
&& \mbox{} +  2g'^{kl} g^{ij} ( \nabla_i \J{k}{r} + \nabla_k \J{i}{r}) \J{j}{q} (\nabla_r g'_{ql}) \\ \nonumber
&& \mbox{} + 2g^{ij} ( \nabla_i \J{q}{r} + \nabla_q \J{i}{r}) (\nabla_r \J{j}{q}) \\ \nonumber
&& \mbox{} + g'^{kl} (\nabla_k \J{l}{r} ) g^{ij} (2\J{j}{s} \nabla_i g'_{rs} + \J{r}{s} \nabla_s g'_{ij})  \\ \nonumber
&& \mbox{} + \J{r}{k} g^{ij} \nabla_i \nabla_j \J{k}{r} - 2 g^{ij} g'^{kl} g'_{qj} \J{r}{q} \nabla_i \nabla_l \J{k}{r} \\ \nonumber
&& \mbox{} - \J{r}{q} g'_{qj} g'^{kl} g^{ij} \nabla_k \nabla_l \J{i}{r} - 2R \\ \label{eqnlaplog}
&& \left. \mbox{} + 2g'^{kl} g'_{ij} g^{pj} R^i_{\, lpk}  - \frac{| \nabla' \tr{g}{g'}|^2}{\tr{g}{g'}} \right),
\end{eqnarray}
where, by definition,
$$|\nabla' \tr{g}{g'}|^2 = g'^{kl} ( \dd{k} \tr{g}{g'}) ( \dd{l} \tr{g}{g'} ),$$
and where we are making use of the equations
\begin{eqnarray*}
g^{ij} g^{kl} \dd{i} \dd{j} g_{kl} + 2 \J{r}{k} g^{ij} \dd{k} \dd{j} \J{i}{r}
& = & - 2R 
\end{eqnarray*}
and
\begin{eqnarray*}
2 J'_{jr} g'^{kl} g^{ij} \dd{i} \dd{l} \J{k}{r} + g'^{kl} g'_{ij} \dd{k}
\dd{l} g^{ij} & = & - 2g^{ij} g'^{kl} g'_{qj} \J{r}{q} \nabla_i \nabla_l
\J{k}{r} +2g'^{kl} g'_{ij} g^{pj} R^i_{\ lpk},
\end{eqnarray*}
which both hold only at a point.

Notice now that (\ref{eqnlaplog}) is an equation of tensors.  Since we are going
to apply the maximum principle we need to obtain a good lower
bound for the right hand side of (\ref{eqnlaplog}).  
We have to deal with
the bad terms that involve derivatives of $g'$ and are not nonnegative:  namely, the third, fifth and last terms. First, we need a lemma.

\begin{lemma} \label{lemmaalphabeta}
Define tensors $\alpha_{ijp}$ and $\beta_{ijp}$ by
\begin{eqnarray*}
\alpha_{ijp} & = & J'_{pl} (\nabla_i \J{j}{l} - \nabla_j \J{i}{l}) + (\nabla_p \J{l}{k}) ( g'_{kj} \J{i}{l} - g'_{ki} \J{j}{l}) \\ 
&& \mbox{} + g'_{kp} ( (\nabla_l \J{j}{k}) \J{i}{l} - (\nabla_l \J{i}{k}) \J{j}{l}), \\ 
\beta_{ijp} & = & g'_{kl}((\nabla_j \J{i}{l}) \J{p}{k} - (\nabla_j \J{p}{l})
\J{i}{k}) - g'_{kj} ( (\nabla_p \J{l}{k}) \J{i}{l} - (\nabla_i \J{l}{k})\J{p}{l})\\
&& \mbox{} - (\nabla_l \J{j}{k}) ( g'_{kp} \J{i}{l} - g'_{ki} \J{p}{l}). 
\end{eqnarray*}
Then
\begin{enumerate}
\item[(i)] $\displaystyle{2 \PP{r}{s}{i}{j} \nabla_r g'_{sp} - 2 \PP{r}{s}{j}{i} \nabla_r g'_{sp} = \alpha_{ijp}}$
\item[(ii)] 
$\displaystyle{2\QQ{r}{s}{i}{j} \nabla_r g'_{sp} - 2\QQ{r}{s}{p}{j} \nabla_r g'_{si} = \beta_{ijp}}.$
\end{enumerate}
\end{lemma}
Notice that in the K\"ahler case the tensors $\alpha$ and $\beta$ vanish identically and the lemma states that $\PP{r}{s}{i}{j} \nabla_r g'_{sp}$ is symmetric in $i$ and $j$ while 
$\QQ{r}{s}{i}{j}\nabla_r g'_{sp}$ is symmetric in $i$ and $p$.  It is important here that $\alpha$ and $\beta$ do not contain derivatives of $g'$.

\begin{proof}
We prove (i).  The proof of (ii) is similar.
From
$$\nabla_i J'_{jk} + \nabla_j J'_{ki} + \nabla_k J'_{ij}=0,$$
we have
$$\J{j}{q} \nabla_i g'_{qk} + \J{k}{q} \nabla_j g'_{qi} + \J{i}{q} \nabla_k g'_{qj} +g'_{qk} \nabla_i \J{j}{q} + g'_{qi} \nabla_j \J{k}{q} + g'_{qj}
 \nabla_k \J{i}{q}=0.$$
Multiplying by $\J{p}{j}$ we obtain
\begin{eqnarray*}
 \nabla_i g'_{pk} & = & \J{p}{j} \J{k}{q} \nabla_j g'_{qi}
+ \J{p}{j} \J{i}{q} \nabla_k g'_{jq} + g'_{qk} \J{p}{j} \nabla_i \J{j}{q}\\
&& \mbox{} + g'_{qi} \J{p}{j} \nabla_j \J{k}{q} + g'_{qj} \J{p}{j} \nabla_k
\J{i}{q} \\
& = & - 2 \PP{r}{s}{p}{k} \nabla_r g'_{si} + \nabla_p g'_{ki} + \nabla_k g'_{pi}
- g'_{qj} \J{i}{q} \nabla_k \J{p}{j}\\
&& \mbox{}  + g'_{qk} \J{p}{j} \nabla_i
\J{j}{q} + g'_{qi} \J{p}{j} \nabla_j \J{k}{q},
\end{eqnarray*}
where we have made use of the identity $\PP{j}{q}{p}{i} g'_{jq}
=0$.
Then (i) follows easily.  Q.E.D.
\end{proof}

We return to equation (\ref{eqnlaplog}).
The third term on the right hand side (ignoring the factor of $1/\tr{g}{g'})$
 can be written
\begin{eqnarray} \nonumber
\lefteqn{2 g'^{kl} g^{ij} (\nabla_i \J{k}{r} + \nabla_k \J{i}{r}) \J{j}{q} \nabla_r g'_{ql}} \\ \nonumber
& = & 2 g'^{kl} g^{ij} g^{rs} (\nabla_i J_{kr} + \nabla_k J_{ir}) \nabla_s
J'_{jl} \\ && \mbox{} - 2g'^{kl} g^{ij} g^{rs} (\nabla_i J_{kr} + \nabla_k J_{ir} ) g'_{ql}
\nabla_s \J{j}{q}. \label{eqncalc2}
\end{eqnarray}
Making use of the identity $\QQ{a}{b}{i}{r} \nabla_k J_{ab} =0$ we rewrite
the first term on the right hand side of (\ref{eqncalc2}) as
\begin{eqnarray} \nonumber
\lefteqn{ 2 g'^{kl} g^{ij} g^{rs} (\nabla_i J_{kr} + \nabla_k J_{ir}) \nabla_s
J'_{jl}} \\ \nonumber
& = & 4 g'^{kl} g^{ij} g^{rs} ( \nabla_k J_{ir} ) (\nabla_s J'_{jl}) + 2
g'^{kl} g^{ij} g^{rs} (\nabla_r J_{ki}) ( \nabla_s J'_{jl}) \\ \nonumber
& = & 4 g'^{kl} g^{ij} g^{rs} \PP{a}{b}{i}{r} (\nabla_k J_{ab}) ( \nabla_s
J'_{jl}) + 2g'^{kl} g^{ij} g^{rs} \PP{a}{b}{k}{i} (\nabla_r J_{ab}) (\nabla_s
J'_{jl}) \\ \label{eqnIplusII}
& = & (I) + (II),
\end{eqnarray}
where 
$$(I) = 4 g'^{kl} g^{ij} g^{rs} (\nabla_k J_{ir}) \PP{a}{b}{j}{s} (\nabla_b
J'_{al})$$
and
$$(II) = 2 g'^{kl} g^{ij} g^{rs} (\nabla_r J_{ik}) \PP{a}{b}{j}{l} (\nabla_s
J'_{ab}).$$
To deal with $(I)$, first note that
\begin{eqnarray} \nonumber
\PP{a}{b}{j}{s} \nabla_b J'_{al} & = & \PP{a}{b}{s}{j} \nabla_a J'_{bl} \\
\nonumber
& = & - \PP{a}{b}{s}{j} (\nabla_a g'_{bq}) \J{l}{q} - \PP{a}{b}{s}{j} g'_{bq}
\nabla_a \J{l}{q} \\ \nonumber
& = & - \PP{a}{b}{j}{s} (\nabla_a g'_{bq}) \J{l}{q} - \frac{1}{2} \alpha_{sjq} \J{l}{q}
- \PP{a}{b}{s}{j} g'_{bq} \nabla_a \J{l}{q} \\ \nonumber
& = & \PP{a}{b}{s}{j} \nabla_b J'_{al} + \PP{a}{b}{j}{s} g'_{bq} \nabla_a
\J{l}{q} - \frac{1}{2} \alpha_{sjq} \J{l}{q}
 - \PP{a}{b}{s}{j} g'_{bq} \nabla_a \J{l}{q}.
\end{eqnarray}
That is, the term $\PP{a}{b}{j}{s} \nabla_b J'_{al}$ is symmetric in $j$
and $s$ modulo terms that don't involve derivatives of $g'$.  Then
\begin{eqnarray} \nonumber
(I) & = & - 4 g'^{kl} g^{ij} g^{rs} ( \nabla_k J_{ri}) \PP{a}{b}{s}{j} \nabla_b
J'_{al} \\ \nonumber
&& \mbox{} + 4 g'^{kl} g^{ij} g^{rs} (\nabla_k J_{ir}) (\PP{a}{b}{j}{s} g'_{bq}
\nabla_a \J{l}{q} - \frac{1}{2}\alpha_{sjq} \J{l}{q} - \PP{a}{b}{s}{j} g'_{bq} \nabla_a
\J{l}{q}),
\end{eqnarray}
and hence, interchanging the indices $i$ and $j$ with $r$ and $s$ in the
first term, we have
\begin{equation} \label{eqnI}
(I) = 2g'^{kl} g^{ij} g^{rs} (\nabla_k J_{ir}) ( \PP{a}{b}{j}{s} g'_{bq} \nabla_a
\J{l}{q} - \frac{1}{2} \alpha_{sjq} \J{l}{q} - \PP{a}{b}{s}{j} g'_{bq} \nabla_a \J{l}{q}).
\end{equation}
For $(II)$, calculate
\begin{eqnarray} \nonumber
\PP{a}{b}{j}{l} \nabla_s J'_{ab} & = & \frac{1}{2} \left( \nabla_s J'_{jl} -
\J{j}{a} \J{l}{b}\nabla_s J'_{ab} \right) \\ \nonumber
& = & \frac{1}{2} \left( \nabla_s J'_{jl} -  \nabla_s J'_{jl} + J'_{ab}
\J{l}{b} \nabla_s \J{j}{a} + J'_{ab} \J{j}{a} \nabla_s \J{l}{b} \right) \\ \label{eqnII}
& = & \frac{1}{2} J'_{ab} \left( \J{l}{b} \nabla_s \J{j}{a} + \J{j}{a} \nabla_s
\J{l}{b} \right). 
\end{eqnarray}
Then from equations (\ref{eqncalc2}), (\ref{eqnIplusII}), (\ref{eqnI}) and
(\ref{eqnII}) we have the following expression for the third term of (\ref{eqnlaplog}):
\begin{eqnarray} \nonumber
\lefteqn{2 g'^{kl} g^{ij} (\nabla_i \J{k}{r} + \nabla_k \J{i}{r}) \J{j}{q} \nabla_r g'_{ql}} \\ \nonumber
& = & 2g'^{kl} g^{ij} g^{rs} (\nabla_k J_{ir}) ( \PP{a}{b}{j}{s} g'_{bq}
\nabla_a \J{l}{q} - \frac{1}{2}\alpha_{sjq} \J{l}{q} - \PP{a}{b}{s}{j} g'_{bq} \nabla_a
\J{l}{q}) \\ \nonumber
&& \mbox{} + g'^{kl} g^{ij} g^{rs} (\nabla_r J_{ik}) J'_{ab} (\J{l}{b} \nabla_s
\J{j}{a} + \J{j}{a} \nabla_s \J{l}{b}) \\ \label{eqn3rdterm}
&& \mbox{} - 2g'^{kl} g^{ij} g^{rs} ( \nabla_i
J_{kr} + \nabla_k J_{ir}) g'_{ql} \nabla_s \J{j}{q}.
\end{eqnarray}

We deal now
with the fifth term on the right hand side of 
(\ref{eqnlaplog}).  Calculate:
\begin{eqnarray} \nonumber
\lefteqn{g'^{kl} (\nabla_k \J{l}{r} ) g^{ij} (2\J{j}{s} \nabla_i g'_{rs} + \J{r}{s} \nabla_s g'_{ij})} \\ \nonumber
& = & g'^{kl} (\nabla_k \J{l}{r}) g^{ij} ( 2 \nabla_i J'_{jr} - \J{r}{s}
\nabla_s (J'_{pj} \J{i}{p})) \\ \nonumber
& = & g'^{kl} (\nabla_k \J{l}{r}) g^{ij} ( 2 \nabla_i J'_{jr} - \J{r}{s}
\J{i}{p} \nabla_p J'_{sj} - \J{r}{s} \J{i}{p} \nabla_j J'_{ps} - \J{r}{s} J'_{pj}
\nabla_s \J{i}{p}) \\ \nonumber
& = & g'^{kl} (\nabla_k \J{l}{r}) g^{ij} (2 \nabla_i J'_{jr} + \J{i}{a} \J{r}{b}
\nabla_a J'_{jb} -\nabla_j J'_{ir}
 + \J{i}{p} J'_{ps} \nabla_j \J{r}{s} \\ \nonumber
 && \mbox{}
  - \J{r}{s} J'_{pj} \nabla_s
\J{i}{p}) \\ \nonumber
& = & g'^{kl} (\nabla_k \J{l}{r}) g^{ij} ( 2 \QQ{a}{b}{i}{r} \nabla_a J'_{jb}
+ \J{i}{p} J'_{ps} \nabla_j \J{r}{s} - \J{r}{s} J'_{pj} \nabla_s \J{i}{p})
\\ \nonumber
& = & g'^{kl} (\nabla_k \J{l}{r})g^{ij} ( 2 \QQ{a}{b}{i}{r} (\nabla_a g'_{bq})
\J{j}{q} + 2 \QQ{a}{b}{i}{r} g'_{bq} \nabla_a \J{j}{q} - g'_{is} \nabla_j
\J{r}{s} \\  \nonumber
&& \mbox{} - \J{r}{s} J'_{pj} \nabla_s \J{i}{p} ) \\ \nonumber
& = & g'^{kl} ( \nabla_k \J{l}{r}) ( \frac{1}{2} \beta_{irq} J^{iq} + 2g^{ij} \QQ{a}{b}{i}{r}
g'_{bq} \nabla_a \J{j}{q} 
 - g^{ij} g'_{is} \nabla_j \J{r}{s} \\ && \mbox{} \label{eqn5thterm}
 - g^{ij} \J{r}{s} J'_{pj} \nabla_s
\J{i}{p}), 
\end{eqnarray}
where, for the last line, we have used the identity
$$2 \QQ{a}{b}{i}{r} (\nabla_a g'_{bq}) J^{iq} = \frac{1}{2}\beta_{irq} J^{iq},$$
which follows immediately from Lemma \ref{lemmaalphabeta}.

Finally, we must deal with the last term on the right hand side of (\ref{eqnlaplog}). We do this by making use of the good second term.  

\begin{lemma} \label{lemmatrgg'} There exists a constant $C'$ depending only on $\sup_M |F|$ and $\| \emph{Rm}(g) \|_{C^0(g)}$ such that
\begin{eqnarray} \nonumber
\frac{|  \nabla' (\emph{tr}_g{g'})|^2}{\emph{tr}_g{g'}} & \le &  g'^{kl} g'^{ij} g^{pq} (\nabla_p g'_{ik}) ( \nabla_q g'_{jl}) + C' (\emph{tr}_{g}{g'}) (\emph{tr}_{g'}{g}).  \end{eqnarray}
\end{lemma}

To prove this lemma we will work in a coordinate system $(x^1, \ldots, x^4)$ centred at a point $p$ such that the first derivatives of the metric $g$ vanish at $p$ and 
$$\J{1}{3} = \J{2}{4}=1= -\J{3}{1}= -\J{4}{2},$$ and all other entries of the matrix $(\J{i}{j})$ are zero at $p$.
Define local vector fields 
$$Z_{\alpha} = \frac{1}{2} \left( \dd{\alpha} - \sqrt{-1} \J{\al}{i} \dd{i} \right),$$
for $\alpha=1,2$.  Set
$$\G{\al}{\beta} = g(Z_{\al}, \overline{Z}_{\beta})$$ and 
$$\Gp{\al}{\be} = g'(Z_{\al}, \overline{Z}_{\beta}),$$
for $\al, \beta =1,2$.
Then we make a linear change in the coordinates $(x^1, \ldots, x^4)$ so that, in addition to the above conditions at $p$, we also impose that  $\G{\al}{\beta} = \delta_{\al \be}$ and that $\Gp{\al}{\be}$ be diagonal.  
  Notice that in these coordinates, the first derivatives of the $g_{ij}$ vanish at $p$, but in general, the first derivatives of the $\G{\al}{\beta}$ will not.  

It will also be useful to consider the local vector fields
$$W_{\alpha} = \frac{1}{2} \left( - \J{A}{i} \dd{i} - \sqrt{-1} \dd{A} \right),$$
for $\alpha = 1,2$, where we are setting $A= \alpha+2$.  Notice that, at $p$, $W_{\alpha} = Z_{\alpha}$.
In the sequel, we will use the indices  $A, B, C,D, M,N $ to denote $\alpha+2, \beta+2, \gamma + 2, \delta+2, \mu+2, \nu+2$.
Set
$$\GT{\al}{\beta} = g(W_{\al}, \overline{W}_{\beta})$$ and
$$\GTp{\al}{\beta} = g'(W_{\al}, \overline{W}_{\beta}).$$
Observe that
$$\G{\al}{\beta} = \frac{1}{2} \left( g_{\al \be} - \sqrt{-1} J_{\al \beta} \right), \qquad 
\Gp{\al}{\beta} = \frac{1}{2} \left( g'_{\al \be} - \sqrt{-1} J'_{\al \beta} \right),$$
and
$$\GT{\al}{\beta} = \frac{1}{2} \left( g_{AB} - \sqrt{-1} J_{AB} \right), \qquad 
\GTp{\al}{\beta} = \frac{1}{2} \left( g'_{AB} - \sqrt{-1} J'_{AB} \right).$$
At the point $p$, $G = \tilde{G}$ and $G' = \tilde{G}'$; this fact will
be used later in the proof of Lemma \ref{lemmatrgg'}.
Notice also that 
$g_{ij} = 2\delta_{ij}$ and that $g'_{ij}$ is diagonal.
A final word about notation: when we are using the local vector fields $Z_{\al}$ and $\overline{Z}_{\al}$ as differential operators, we will instead write $D_{\al}$ and $D_{\ov{\al}}$ respectively.  We require some preliminary results before we prove Lemma \ref{lemmatrgg'}.

\begin{lemma} \label{lemmaDG}
\begin{eqnarray} \label{eqnDG}
D_{\gamma} \Gp{\al}{\be} & = & D_{\al} \Gp{\ga}{\be} + \frac{1}{4} \left( a_{\ga \al \beta} + \sqrt{-1} b_{\ga \al \be}\right) \\ \label{eqnDGT}
\textrm{and} \quad D_{\gamma} \GTp{\al}{\be} & = & D_{\al} \GTp{\ga}{\be} + \frac{1}{4} \left( a_{\ga A B} + \sqrt{-1} b_{\ga A B}\right), 
\end{eqnarray}
where
\begin{eqnarray*}
a_{\ga i j} & = & ( \dd{i} \J{\ga}{k} - \dd{\ga} \J{i}{k} ) J'_{k j} + (\dd{j} \J{\ga}{k}) J'_{i k} - (\dd{j} \J{i}{k} ) J'_{\gamma k} \\
\textrm{and} \quad b_{\ga i j} & = & g'_{k i} \dd{j} \J{\ga}{k} - g'_{k \ga} \dd{j} \J{i}{k} + J'_{j l} ( \J{i}{k} \dd{k} \J{\ga}{l} - \J{\ga}{k} \dd{k} \J{i}{l}).
\end{eqnarray*}
\end{lemma}
\begin{proof} We will just prove (\ref{eqnDG}), since the proof of (\ref{eqnDGT}) is similar. Calculate
\begin{eqnarray} \nonumber
D_{\ga} \Gp{\al}{\be} & = & \frac{1}{4} ( \dd{\ga} - \sqrt{-1} \J{\ga}{i} \dd{i}) (g'_{\al \be} - \sqrt{-1} J'_{\al \be}) \\ \nonumber
& = & \frac{1}{4} ( \dd{\ga} g'_{\al \be} - \J{\ga}{i} \dd{i} J'_{\al \be} - \sqrt{-1} \J{\ga}{i} \dd{i} g'_{\al \be} - \sqrt{-1} \dd{\ga} J'_{\al \be} ) \\ \nonumber
& = & \frac{1}{4} ( \dd{\ga} g'_{\al \be} - \J{\ga}{i} \dd{\al} J'_{i\beta} - \J{\ga}{i} \dd{\beta} J'_{\al i} - \sqrt{-1} \J{\ga}{i} \dd{i} g'_{\al \be} \\ \nonumber
&& \mbox{} - \sqrt{-1} \dd{\al} J'_{\ga \beta} - \sqrt{-1} \dd{\be} J'_{\al \ga}) \\ \nonumber
&  = & \frac{1}{4} (  \dd{\ga} g'_{\al \beta} + \dd{\al} g'_{\ga \beta} + (\dd{\al} \J{\ga}{i}) J'_{i \be} - \J{\al}{i} \dd{\be} J'_{\ga i} 
 - (\dd{\be} \J{\al}{i}) J'_{\ga i} \\ \nonumber && \mbox{} + (\dd{\be} \J{\ga}{i} ) J'_{\al i} - \sqrt{-1} \J{\ga}{i} \dd{i} (J'_{\be p} \J{\al}{p}) - \sqrt{-1} \dd{\al} J'_{\ga \be} - \sqrt{-1} \dd{\be} J'_{\al \ga}) \\ \nonumber
& = & \frac{1}{4} ( \dd{\ga} g'_{\al \be} + \dd{\al} g'_{\ga \be} + (\dd{\al} \J{\ga}{i} ) J'_{i\be} - \J{\al}{i} \dd{i} J'_{\ga \be} - \J{\al}{i} \dd{\ga} J'_{\beta i}  - (\dd{\be} \J{\al}{i} ) J'_{\ga i} \\ \nonumber
&& \mbox{}+ (\dd{\be} \J{\ga}{i}) J'_{\al i} - \sqrt{-1} \J{\ga}{i} \J{\al}{p} \dd{i} J'_{\be p} - \sqrt{-1} \J{\ga}{i} J'_{\beta p } \dd{i} \J{\al}{p} \\ \label{eqncalcDGp}
&& \mbox{} - \sqrt{-1} \dd{\al} J'_{\ga \be} - \sqrt{-1} \dd{\be} J'_{\al \ga}), 
\end{eqnarray}
where to go from the third to the fourth lines we have used the simple identity
$$\J{\ga}{i} \dd{\be} J'_{\al i} + (\dd{\be} \J{\ga}{i}) J'_{\al i}
 = (\dd{\be} \J{\al}{i} ) J'_{\ga i} + \J{\al}{i} \dd{\be} J'_{\ga i} .$$
Notice that
\begin{eqnarray*}
\J{\ga}{i} \J{\al}{p} \dd{i} J'_{\be p} & = & \J{\ga}{i} \J{\al}{p} \dd{\be} J'_{ip} + \J{\ga}{i} \J{\al}{p} \dd{p} J'_{\be i} \\
& = & \dd{\be} J'_{\ga \al} - J'_{ip} \J{\al}{p} \dd{\be} \J{\ga}{i} - J'_{ip} \J{\ga}{i} \dd{\be} \J{\al}{p} + \J{\al}{i} \J{\ga}{p} \dd{i} J'_{\be p}.
\end{eqnarray*}
In (\ref{eqncalcDGp}) this gives us
\begin{eqnarray*}
D_{\ga} \Gp{\al}{\be} & = & \frac{1}{4} (\dd{\ga} g'_{\al \be} + \dd{\al} g'_{\ga \be} + (\dd{\al} \J{\ga}{i}) J'_{i\be} - \J{\al}{i} \dd{i} J'_{\ga \be} - \J{\al}{i} \dd{\ga} J'_{\be i} \\
&& \mbox{} - (\dd{\be} \J{\al}{i}) J'_{\ga i} + (\dd{\be} \J{\ga}{i} ) J'_{\al i} + \sqrt{-1} J'_{ip} \J{\al}{p} \dd{\be} \J{\ga}{i} + \sqrt{-1} J'_{ip} \J{\ga}{i} \dd{\be} \J{\al}{p} \\
&& \mbox{} - \sqrt{-1} \J{\al}{i} \J{\ga}{p} \dd{i} J'_{\beta p} - \sqrt{-1} \J{\ga}{i} J'_{\be p} \dd{i} \J{\al}{p} - \sqrt{-1} \dd{\al} J'_{\ga \be}) \\
& = & \frac{1}{4} ( \dd{\al} g'_{\ga \be} + (\dd{\al} \J{\ga}{i} ) J'_{i \be} - \J{\al}{i} \dd{i} J'_{\ga \be} + (\dd{\ga} \J{\al}{i}) J'_{\be i} - (\dd{\be} \J{\al}{i} ) J'_{\ga i} \\ && \mbox{} + (\dd{\be} \J{\ga}{i}) J'_{\al i} + \sqrt{-1} g'_{i\al} \dd{\be} \J{\ga}{i} - \sqrt{-1} g'_{p \ga} \dd{\be} \J{\al}{p} 
 - \sqrt{-1} \J{\al}{i} \dd{i} g'_{\be \ga} \\ && \mbox{} +\sqrt{-1} \J{\al}{i} J'_{\be p} (\dd{i} \J{\ga}{p}) - \sqrt{-1} \J{\ga}{i} J'_{\be p} (\dd{i} \J{\al}{p}) - \sqrt{-1} \dd{\al} J'_{\ga \be}) \\
 & = & D_{\al} \Gp{\ga}{\be} + \frac{1}{4} \left( (\dd{\al} \J{\ga}{i} - \dd{\ga} \J{\al}{i}) J'_{i \be} - (\dd{\be} \J{\al}{i}) J'_{\ga i} + (\dd{\be} \J{\ga}{i}) J'_{\al i} \right. \\
 && \mbox{} + \left. \sqrt{-1} ( g'_{i \al} \dd{\be} \J{\ga}{i} - g'_{i \ga} \dd{\be} \J{\al}{i}) + \sqrt{-1} J'_{\be p} (\J{\al}{i} \dd{i} \J{\ga}{p} - \J{\ga}{i} \dd{i} \J{\al}{p}) \right),
\end{eqnarray*}
as required.  Q.E.D.
\end{proof}

We will also need the following result.

\begin{lemma} \label{lemmaDtrGp} At the point $p$ we have
\begin{eqnarray}  \label{eqnDtrGp}
2 D_{\gamma} (\GI{\al}{\be} \Gp{\al}{\be}) & = & 2 \GI{\al}{\be} D_{\ga} \Gp{\al}{\be} = D_{\gamma} ( g^{ij} g'_{ij}) + E_{\ga} \quad \textrm{and}\\
2 D_{\gamma} (\GTI{\al}{\be} \GTp{\al}{\be}) & = & 2 \GTI{\al}{\be} D_{\ga} \GTp{\al}{\be} = D_{\gamma} (g^{ij} g'_{ij}) - E_{\ga}, 
\end{eqnarray}
where $\displaystyle{E_{\gamma} = \sum_{i=1}^2 g'_{rs} (D_{\gamma} \J{i}{r}) \J{i}{s}}$.
\end{lemma}

We are using here the usual notation for the inverse metrics $G^{-1}$ and $\tilde{G}^{-1}$.  Note also that repeated greek indices $\alpha, \beta, \ldots$ are used to denote a sum from 1 to 2, whereas repeated lower case roman letters $i, j, \ldots$ denote the usual sum from 1 to 4, unless otherwise indicated.

\begin{proof}
We will prove just (\ref{eqnDtrGp}).  Calculate
\begin{eqnarray} \nonumber
2D_{\gamma} ( \GI{\al}{\be} \Gp{\al}{\be})  & = & - \frac{1}{2} \GI{\al}{\mu} \GI{\nu}{\be}( - \J{\gamma}{i} \dd{i} J_{\nu \mu} - \sqrt{-1} \dd{\gamma} J_{\nu \mu}) \Gp{\al}{\be} + 2\GI{\al}{\be} D_{\gamma} \Gp{\al}{\be} \\ \nonumber
& = & 2\GI{\al}{\be} D_{\gamma} \Gp{\al}{\be} \\ \label{eqncalcDtrGp1}
& = & \frac{1}{2} \GI{\al}{\be} (\dd{\gamma} g'_{\al \be} - \sqrt{-1} \J{\ga}{k} \dd{k} g'_{\al \be} ).
\end{eqnarray}
But
\begin{eqnarray} \nonumber
D_{\gamma} ( g^{ij} g'_{ij} ) & = & \frac{1}{2} g^{ij} ( \dd{\gamma} - \sqrt{-1} \J{\ga}{k} \dd{k}) g'_{ij} \\ \nonumber
& = & \frac{1}{4} \sum_{i=1}^4 \left( \dd{\gamma} g'_{ii} - \sqrt{-1} \J{\gamma}{k} \dd{k} g'_{ii} \right) \\ \nonumber
& = & \frac{1}{2} \sum_{i=1}^2 \left( \dd{\gamma} g'_{ii} - g'_{rs} (\dd{\gamma} \J{i}{r}) \J{i}{s} - \sqrt{-1} \J{\gamma}{k} \dd{k} g'_{ii} \right. \\ \label{eqncalcDtrGp2}
&& \mbox{} \left.+ \sqrt{-1} \J{\gamma}{k} g'_{rs} (\dd{k} \J{i}{r}) \J{i}{s} \right),
\end{eqnarray}
where we have used the fact that, at $p$,
$$\dd{k} g'_{AB} - \dd{k} g'_{\alpha \beta} = - g'_{rs} (\dd{k} \J{\al}{r}) \J{\be}{s} - g'_{rs} (\dd{k} \J{\be}{s}) \J{\al}{r}.$$
Comparing (\ref{eqncalcDtrGp1}) and (\ref{eqncalcDtrGp2}) gives (\ref{eqnDtrGp}).
Q.E.D.
\end{proof}

We need one more lemma before we can prove Lemma \ref{lemmatrgg'}.

\begin{lemma} \label{lemmaab}
At $p$,
$$g^{ij} a_{\gamma ij} = 0 \quad \textrm{and} \quad g^{ij} b_{\gamma i j}=0.$$
\end{lemma}
\begin{proof}
For the first equation, calculate at $p$:
\begin{eqnarray*}
g^{ij} a_{\gamma ij} & = & g^{ij} ( \nabla_{i} \J{\gamma}{k}) J'_{kj} - g^{ij} (\nabla_{\gamma} \J{i}{k}) J'_{kj} \\ && \mbox{} + g^{ij} (\nabla_{j} \J{\gamma}{k}) J'_{ik} - g^{ij} (\nabla_{j} \J{i}{k}) J'_{\gamma k} \\
& = & -g^{ij} (\nabla_{\gamma} \J{i}{k}) J'_{kj} \\
& = & -g^{ij} (\nabla_{\gamma} \J{i}{k}) \J{k}{s} g'_{sj} \\
& = & g^{ij} \J{i}{k} (\nabla_{\gamma} \J{k}{s}) g'_{sj} \\
& = & - g^{ik} \J{i}{j} (\nabla_{\gamma} \J{k}{s}) g'_{sj} \\
& = & 0, 
\end{eqnarray*}
by symmetry.  For the second, calculate:
\begin{eqnarray*}
g^{ij} b_{\gamma ij} & = & g^{ij} g'_{ki} (\nabla_j \J{\ga}{k}) - g^{ij} g'_{k \ga} (\nabla_j \J{i}{k}) \\ && \mbox{}+ g^{ij} J'_{jl} \J{i}{k} (\nabla_k \J{\ga}{l})  - g^{ij} J'_{jl} \J{\ga}{k} (\nabla_k \J{i}{l}) \\
& = & g^{ij} J'_{kl} \J{i}{l} (\nabla_j \J{\ga}{k}) - g^{ik} J'_{jl} \J{i}{j} (\nabla_k \J{\ga}{l}) \\
&& \mbox{} - g^{ij} J'_{jl} \J{\ga}{k} (\nabla_i \J{k}{l}) - g^{ij} J'_{jl} \J{\ga}{k} (\nabla^l J_{ik}) \\
& = & 2 g^{ij} J'_{kl} \J{i}{l} (\nabla_j \J{\ga}{k}) - 2 \J{\ga}{k} J'_{jl} (\nabla^j \J{k}{l}).
\end{eqnarray*}
Notice that the second term vanishes since
$$J'_{jl} \nabla^j \J{k}{l} = \frac{1}{2} g'_{lp} \J{j}{p} \nabla_k J^{jl}=0.$$
But also the first term can be written
\begin{eqnarray*}
g^{ij} J'_{kl} \J{i}{l} (\nabla_j \J{\ga}{k}) & = & g^{ij} \J{k}{p} g'_{pl} \J{i}{l} (\nabla_j \J{\ga}{k})  =  - J'_{ip} \J{\ga}{k} (\nabla^i \J{k}{p})  =  0,
\end{eqnarray*}
finishing the proof of the lemma.  Q.E.D.
\end{proof}

\bigskip
\noindent
{\bf Proof of Lemma \ref{lemmatrgg'} } \ Using Lemma \ref{lemmaDtrGp}  we have
\begin{eqnarray*}
\lefteqn{g'^{kl} \dd{k} (g^{ij} g'_{ij}) \dd{l} (g^{pq} g'_{pq}) } \\& = & 2 \GpI{\ga}{\de} D_{\ga} (g^{ij} g'_{ij}) D_{\ov{\de}} (g^{pq}g'_{pq}) \\ 
& = & \GpI{\ga}{\de}( 2 D_{\ga}( \GI{\al}{\be} \Gp{\al}{\be}) - E_{\ga})(2 D_{\ov{\de}} (\GI{\mu}{\nu} \Gp{\mu}{\nu}) - \ov{E}_{\de}) \\
&& \mbox{} + \GpI{\ga}{\de}( 2 D_{\ga}( \GTI{\al}{\be} \GTp{\al}{\be}) + E_{\ga})(2 D_{\ov{\de}} (\GTI{\mu}{\nu} \GTp{\mu}{\nu}) + \ov{E}_{\de}) \\
& = & 4 \GpI{\ga}{\de}( \GI{\al}{\be} D_{\ga} \Gp{\al}{\be})( \GI{\mu}{\nu} D_{\ov{\de}} \Gp{\mu}{\nu}) \\
&& \mbox{} + 4 \GpI{\ga}{\de}( \GTI{\al}{\be} D_{\ga} \GTp{\al}{\be})( \GTI{\mu}{\nu} D_{\ov{\de}} \GTp{\mu}{\nu}) - 2 \GpI{\ga}{\de} E_{\ga} \ov{E}_{\de}.
\end{eqnarray*}
From Lemma \ref{lemmaDG},
\begin{eqnarray*}
\lefteqn{g'^{kl} \dd{k} (g^{ij} g'_{ij}) \dd{l} (g^{pq} g'_{pq}) } \\
& = & 4 \GpI{\ga}{\de} (\GI{\al}{\be} D_{\al} \Gp{\ga}{\be})( \GI{\mu}{\nu} D_{\ov{\nu}} \Gp{\mu}{\de})  \\
&& \mbox{} + 4 \GpI{\ga}{\de} (\GI{\al}{\be} D_{\al} \GTp{\ga}{\be})( \GI{\mu}{\nu} D_{\ov{\nu}} \GTp{\mu}{\de}) \\
&& \mbox{} + 2\textrm{Re} \left\{ \GpI{\ga}{\de} \GI{\al}{\be}( a_{\ga \al \be} + \sqrt{-1} b_{\ga \al \be}) \GI{\mu}{\nu} D_{\ov{\nu}} \Gp{\mu}{\de} \right. \\
&& \mbox{}+ \left. \GpI{\ga}{\de} \GI{\al}{\be}( a_{\ga A B} + \sqrt{-1} b_{\ga A B}) \GI{\mu}{\nu} D_{\ov{\nu}} \GTp{\mu}{\de} \right\} \\
&& \mbox{} + \frac{1}{4} \GpI{\ga}{\de} \GI{\al}{\be}(a_{\ga \al \be} + \sqrt{-1} b_{\ga \al \be}) \GI{\mu}{\nu} (a_{\de \mu \nu} - \sqrt{-1} b_{\de \mu \nu}) \\
&& \mbox{} + \frac{1}{4} \GpI{\ga}{\de} \GI{\al}{\be}(a_{\ga AB} + \sqrt{-1} b_{\ga AB}) \GI{\mu}{\nu} (a_{\de MN} - \sqrt{-1} b_{\de MN}) \\
&& \mbox{} - 2 \GpI{\ga}{\de} E_{\ga} \ov{E}_{\de}.
\end{eqnarray*}
But
\begin{eqnarray*}
\GI{\mu}{\nu} D_{\ov{\nu}} \GTp{\mu}{\de} & = & \GI{\mu}{\nu} D_{\ov{\nu}} \Gp{\mu}{\de} + \frac{1}{4}\GI{\mu}{\nu} \left( a_{\de \mu \nu} + a_{\nu MD} - \sqrt{-1} b_{\de \mu \nu} - \sqrt{-1} b_{\nu MD} \right) - \ov{E}_{\de},
\end{eqnarray*}
and so, making use of Lemma \ref{lemmaab},
\begin{eqnarray} \nonumber
\lefteqn{g'^{kl} \dd{k} (g^{ij} g'_{ij}) \dd{l} (g^{pq} g'_{pq}) } \\ \nonumber
& = & 4 \GpI{\ga}{\de} (\GI{\al}{\be} D_{\al} \Gp{\ga}{\be})( \GI{\mu}{\nu} D_{\ov{\nu}} \Gp{\mu}{\de})  \\ \nonumber
&& \mbox{} + 4 \GpI{\ga}{\de} (\GI{\al}{\be} D_{\al} \GTp{\ga}{\be})( \GI{\mu}{\nu} D_{\ov{\nu}} \GTp{\mu}{\de}) \\ \nonumber
&& \mbox{} + \frac{1}{2} \textrm{Re} \left\{ \GpI{\ga}{\de} \GI{\al}{\be}(a_{\ga AB} + \sqrt{-1} b_{\ga AB}) (\GI{\mu}{\nu} (a_{\de \mu \nu } + a_{\nu MD} \right. \\ \nonumber
&& \left. \mbox{} - \sqrt{-1} b_{\de \mu \nu} - \sqrt{-1} b_{\nu MD}) - 4 \ov{E}_{\de} ) \right\} \\ \nonumber
&& \mbox{} + \frac{1}{4} \GpI{\ga}{\de} \GI{\al}{\be}(a_{\ga \al \be} + \sqrt{-1} b_{\ga \al \be}) \GI{\mu}{\nu} (a_{\de \mu \nu} - \sqrt{-1} b_{\de \mu \nu}) \\ \nonumber
&& \mbox{} + \frac{1}{4} \GpI{\ga}{\de} \GI{\al}{\be}(a_{\ga AB} + \sqrt{-1} b_{\ga AB}) \GI{\mu}{\nu} (a_{\de MN} - \sqrt{-1} b_{\de MN}) \\ \label{eqncalclong1}
&& \mbox{} - 2 \GpI{\ga}{\de} E_{\ga} \ov{E}_{\de}.
\end{eqnarray}
We can now apply the Cauchy-Schwartz inequality twice to obtain:
\begin{eqnarray} \nonumber
\lefteqn{\GpI{\ga}{\de} (\GI{\al}{\be} D_{\al} \Gp{\ga}{\be}) (\GI{\mu}{\nu}  D_{\ov{\nu}} \Gp{\mu}{\de})}\\ \nonumber
& = & \sum_{\ga, \al, \mu} \GpI{\ga}{\ga}( D_{\al} \Gp{\ga}{\al}) (D_{\ov{\mu}} \Gp{\mu}{\ga}) \\ \nonumber
& \le & \sum_{\al, \mu} \left( \sum_{\ga} \GpI{\ga}{\ga} | D_{\al} \Gp{\ga}{\al} |^2 \right)^{1/2} \left( \sum_{\ga} \GpI{\ga}{\ga} | D_{\mu} \Gp{\ga}{\mu}|^2 \right)^{1/2} \\ \nonumber
& = & \left( \sum_{\al} \left( \sum_{\ga} \GpI{\ga}{\ga} | D_{\al} \Gp{\ga}{\al}|^2 \right)^{1/2} \right)^2 \\ \nonumber
& = & \left( \sum_{\al} \sqrt{ \Gp{\al}{\al}} \left( \sum_{\ga} \GpI{\ga}{\ga} \GpI{\al}{\al} | D_{\al} \Gp{\ga}{\al}|^2 \right)^{1/2} \right)^2 \\ \nonumber
& \le & \left( \sum_{\al} \Gp{\al}{\al} \right) \sum_{\ga, \al} \GpI{\ga}{\ga} \GpI{\al}{\al} | D_{\al} \Gp{\ga}{\al} |^2 \\ \label{eqncauchyschwartz1}
& \le & 4 (\tr{g}{g'}) \sum_{\ga, \al, \mu} g'^{\ga \ga} g'^{\al \al} g^{\mu \mu} | D_{\mu} \Gp{\ga}{\al} |^2.
\end{eqnarray}
Similarly,
\begin{eqnarray} \nonumber
\lefteqn{\GpI{\ga}{\de} (\GI{\al}{\be} D_{\al} \GTp{\ga}{\be}) (\GI{\mu}{\nu}  D_{\ov{\nu}} \GTp{\mu}{\de})} \\ \label{eqncauchyschwartz2}
& \le & 4 (\tr{g}{g'}) \sum_{\ga, \al, \mu} g'^{\ga \ga} g'^{\al \al} g^{\mu \mu} | D_{\mu} \GTp{\ga}{\al} |^2.
\end{eqnarray}
But
\begin{eqnarray*}
\lefteqn{ 4 \sum_{\ga, \al, \mu} g'^{\ga \ga} g'^{\al \al} g^{\mu \mu} | D_{\mu} \Gp{\ga}{\al}|^2 } \\
& = & \frac{1}{4} \sum_{\ga, \al, \mu} g'^{\ga \ga} g'^{\al \al} g^{\mu \mu} \left\{ ( \dd{\mu} g'_{\ga \al} - \J{\mu}{i} g'_{j \al} (\nabla_i \J{\ga}{j}) - \J{\mu}{i} \J{\ga}{j} (\nabla_i g'_{j \al}))^2 \right. \\
&& \mbox{} \left. + ( \J{\mu}{i} \nabla_i g'_{\ga \al} + g'_{j\al} \nabla_{\mu} \J{\ga}{j} + \J{\ga}{j} \nabla_{\mu} g'_{j \al})^2 \right\} \\
& = & \frac{1}{4} \sum_{\ga, \al, \mu} \left( g'^{\ga \ga} g'^{\al \al} g^{\mu \mu} ( \nabla_{\mu} g'_{\ga \al})^2 + g'^{\ga \ga} g'_{\al \al} g^{MM} ( \nabla_M \J{\ga}{\al})^2 \right. \\
&& \mbox{} + g'^{CC} g'^{\al \al} g^{MM} (\nabla_M g'_{C\al})^2 + g'^{\ga \ga} g'^{\al \al} g^{MM} (\nabla_M g'_{\ga \al})^2 \\
&& \mbox{} + g'^{\ga \ga} g'_{\al \al} g^{\mu \mu} (\nabla_{\mu} \J{\ga}{\al})^2 + g'^{CC} g'^{\al \al} g^{\mu \mu} (\nabla_{\mu} g'_{C\al})^2 \\ 
&& \mbox{} - 2 g'^{\ga \ga} g^{\mu \mu} ( \nabla_{\mu} g'_{\ga \al}) \J{\mu}{i} (\nabla_i \J{\ga}{\al}) - 2 g'^{\ga \ga} g'^{\al \al} g^{\mu \mu} \J{\mu}{i} \J{\ga}{j} (\nabla_{\mu} g'_{\ga \al})(\nabla_i g'_{j\al}) \\
&& \mbox{} + 2 g'^{\ga \ga} \J{\ga}{j} g^{MM} (\nabla_M \J{\ga}{\al}) (\nabla_M g'_{j \al}) + 2 g'^{\ga \ga} g^{\mu \mu} \J{\mu}{i} (\nabla_i g'_{\ga \al} ) (\nabla_{\mu} \J{\ga}{\al}) \\
&& \mbox{} \left. + 2 g'^{\ga \ga} g'^{\al \al} g^{\mu \mu} \J{\mu}{i} \J{\ga}{j} (\nabla_i g'_{\ga \al}) ( \nabla_{\mu} g'_{j \al}) + 2g'^{\ga \ga} g^{\mu \mu} \J{\ga}{j} (\nabla_{\mu} \J{\ga}{\al} )(\nabla_{\mu} g'_{j \al}) \right).
\end{eqnarray*}
A short calculation shows that
$$ \sum_{\ga, \al, \mu} g'^{\ga \ga} g^{\mu \mu} ( \nabla_{\mu} g'_{\ga \al}) \J{\mu}{i} (\nabla_i \J{\ga}{\al}) + g'^{\ga \ga} g'^{\al \al} g^{\mu \mu} \J{\mu}{i} \J{\ga}{j} (\nabla_{\mu} g'_{\ga \al})(\nabla_i g'_{j\al})=0,$$
and similarly that
$$ \sum_{\ga, \al, \mu} g'^{\ga \ga} g^{\mu \mu} \J{\mu}{i} (\nabla_i g'_{\ga \al} ) (\nabla_{\mu} \J{\ga}{\al}) + g'^{\ga \ga} g'^{\al \al} g^{\mu \mu} \J{\mu}{i} \J{\ga}{j} (\nabla_i g'_{\ga \al}) ( \nabla_{\mu} g'_{j \al}) = 0.$$
Hence
\begin{eqnarray} \nonumber
\lefteqn{ 4 \sum_{\ga, \al, \mu} g'^{\ga \ga} g'^{\al \al} g^{\mu \mu} | D_{\mu} \Gp{\ga}{\al}|^2 } \\ \nonumber
 & = & \frac{1}{4} \sum_{\ga, \al=1}^2 \sum_{k=1}^4 \left( g'^{\ga \ga} g'^{\al \al} g^{k k} ( \nabla_{k} g'_{\ga \al})^2 
+ g'^{CC} g'^{\al \al} g^{kk} (\nabla_k g'_{C\al})^2 \right. \\  \label{eqnsumDsquared1}
&& \mbox{} \left.   +  2 g'^{\ga \ga} \J{\ga}{j} g^{kk} (\nabla_k \J{\ga}{\al}) (\nabla_k g'_{j \al}) 
 + g'^{\ga \ga} g'_{\al \al} g^{kk} (\nabla_{k} \J{\ga}{\al})^2 
  \right).
\end{eqnarray}
Similarly, we have
\begin{eqnarray} \nonumber
\lefteqn{ 4 \sum_{\ga, \al, \mu} g'^{\ga \ga} g'^{\al \al} g^{\mu \mu} | D_{\mu} \GTp{\ga}{\al}|^2 } \\ \nonumber
 & = & \frac{1}{4} \sum_{\ga, \al=1}^2 \sum_{k=1}^4 \left( g'^{CC} g'^{AA} g^{k k} ( \nabla_{k} g'_{C A})^2 
+ g'^{\ga \ga} g'^{A A} g^{kk} (\nabla_k g'_{\ga A})^2 \right. \\  \label{eqnsumDsquared2}
&& \mbox{} \left.   -  2 g'^{CC}  g^{kk} (\nabla_k \J{C}{A}) (\nabla_k g'_{\gamma A}) 
 + g'^{C C} g'_{AA} g^{kk} (\nabla_{k} \J{C}{A})^2 
  \right).
\end{eqnarray}
Now observe that, at $p$,
$$\nabla_k \J{C}{A} = - \nabla_k \J{\ga}{\al},$$
and that
$$\nabla_k g'_{\ga A} + \nabla_k g'_{C\al} = J'_{\al i} (\nabla_k \J{C}{i}) + J'_{Ci} (\nabla_k \J{\al}{i}).$$
Using these two simple equalities we obtain at the point $p$, by combining (\ref{eqnsumDsquared1}) and (\ref{eqnsumDsquared2}),
\begin{eqnarray*}
\lefteqn{ 4 \sum_{\ga, \al, \mu} g'^{\ga \ga} g'^{\al \al} g^{\mu \mu} | D_{\mu} \Gp{\ga}{\al}|^2   +  4 \sum_{\ga, \al, \mu} g'^{\ga \ga} g'^{\al \al} g^{\mu \mu} | D_{\mu} \GTp{\ga}{\al}|^2  } \\
& = & \frac{1}{4} g'^{kl} g'^{ij} g^{pq} (\nabla_p g'_{ik}) (\nabla_q g'_{jl}) + \frac{1}{4} \sum g'^{\ga \ga} g'_{\al \al} g^{kk} (\nabla_k \J{\ga}{\al})^2 \\
&& \mbox{} + \frac{1}{4} \sum g'^{CC} g'_{AA} g^{kk} (\nabla_k \J{C}{A})^2 \\
&& \mbox{} - \frac{1}{2} \sum g'^{CC}  g^{kk} (\nabla_k \J{C}{A}) (J'_{\al i} \nabla_k \J{C}{i} + J'_{C i} \nabla_k \J{\al}{i}).
\end{eqnarray*}
Using this, together with (\ref{eqncalclong1}), (\ref{eqncauchyschwartz1}) and (\ref{eqncauchyschwartz2}), we obtain the estimate of Lemma \ref{lemmatrgg'}.
 Q.E.D.

\bigskip

Returning to the proof of the theorem,
we can now combine equations (\ref{eqnlaplog}) with (\ref{eqn3rdterm}) and (\ref{eqn5thterm}) together with Lemma \ref{lemmatrgg'} to obtain
\begin{eqnarray*}
\lefteqn{ \Delta' (\log \tr{g}{g'} ) } \\
 & \ge & \frac{1}{\tr{g}{g'}} \{ 2\Delta F + 2g'^{kl} g^{ij} g^{rs} (\nabla_k J_{ir}) ( \PP{a}{b}{j}{s} g'_{bq}
\nabla_a \J{l}{q} - \frac{1}{2}\alpha_{sjq} \J{l}{q}  \\ && \mbox{} - \PP{a}{b}{s}{j} g'_{bq} \nabla_a
\J{l}{q})  + g'^{kl} g^{ij} g^{rs} (\nabla_r J_{ik}) J'_{ab} (\J{l}{b} \nabla_s
\J{j}{a} + \J{j}{a} \nabla_s \J{l}{b}) \\ 
&& \mbox{} - 2g'^{kl} g^{ij} g^{rs} ( \nabla_i
J_{kr} + \nabla_k J_{ir}) g'_{ql} \nabla_s \J{j}{q} \\
&& \mbox{} + 2g^{ij} ( \nabla_i \J{q}{r} + \nabla_q \J{i}{r}) (\nabla_r \J{j}{q}) \\ 
&& \mbox{} + g'^{kl} ( \nabla_k \J{l}{r}) ( \frac{1}{2} \beta_{irq} J^{iq} + 2g^{ij} \QQ{a}{b}{i}{r}
g'_{bq} \nabla_a \J{j}{q} \\ && \mbox{} 
 - g^{ij} g'_{is} \nabla_j \J{r}{s} 
 - g^{ij} \J{r}{s} J'_{pj} \nabla_s
\J{i}{p})  \\
&& \mbox{} + \J{r}{k} g^{ij} \nabla_i \nabla_j \J{k}{r} - 2 g^{ij} g'^{kl} g'_{qj} \J{r}{q} \nabla_i \nabla_l \J{k}{r} \\ 
&& \mbox{} - \J{r}{q} g'_{qj} g'^{kl} g^{ij} \nabla_k \nabla_l \J{i}{r} - 2R  + 2g'^{kl} g'_{ij} g^{pj} R^i_{\, lpk} \\
&& \mbox{}   - C' (\tr{g}{g'})(\tr{g'}{g}) 
\}.
\end{eqnarray*}
From the Calabi-Yau equation and the arithmetic-geometric means inequality
we have
$$ \tr{g}{g'} \ge 4 \exp{\left( \frac{\inf_M F}{2}\right)} >0.$$ 
Hence, recalling from section 2 that $\|
J \|_{C^2}$ can be bounded in terms of $\| \Rm(g) \|_{C^2}$, we have:
\begin{equation} \label{eqnmax}
\Delta' (\log \tr{g}{g'}) \ge - \frac{A}{2} \tr{g'}{g} - A,
\end{equation}
for a constant $A$ depending only on $\| \Rm(g) \|_{C^2}$, $\sup|F|$ and
the lower bound of $\Delta F$.  We now apply the
maximum principle to the quantity $(\log \tr{g}{g'} -A\phi_1)$.  Suppose that
the maximum is achieved at a point $x_0$ on $M$.  Then at $x_0$ we have
$$\Delta' (\log \tr{g}{g'} - A\phi_1) \le 0.$$
Recall that $\Delta' \phi_1 = 4 - \tr{g'}{g}$.  
At $x_0$ we see that from (\ref{eqnmax})
\begin{eqnarray*}
0 & \ge & \Delta' (\log\tr{g}{g'} -A\phi_1) \\
&\ge & \mbox{} 
 \frac{A}{2} \tr{g'}{g} - 5A,
\end{eqnarray*}
so that $(\tr{g'}{g})(x_0) \le 10$.
On the other hand, from (\ref{eqnCY3}), we have that
$$(\tr{g}{g'})(x_0) \le  10 e^{F(x_0)},$$
and thus at any point $x$ we have
$$\log((\tr{g}{g'})(x)) - A \phi_1(x) \le \log 10 e^{F(x_0)} 
- A\phi_1(x_0).$$
The theorem follows after exponentiating.  Q.E.D.

\bigskip
\bigskip
\pagebreak[3]
\noindent
{\bf 4. H\"older estimate on the metric}
\addtocounter{section}{1}
\setcounter{equation}{0}
\setcounter{theorem}{0}
\setcounter{lemma}{0}
\bigskip

In this section we will prove a H\"older estimate on $g'$ given
a uniform estimate of $g'$, using a modification of the method of Evans \cite{Ev}
and Krylov \cite{Kr} (see also the estimate of Trudinger \cite{Tr2}
and the exposition of Siu \cite{Si}).

\begin{theorem} \label{theoremholder}
Suppose that $g'$ satisfies the equation (\ref{eqnCY2}) and there exists a constant $C_0$
with 
$$C_0^{-1} g \le g' \le C_0 g.$$
Then there exist positive constants $C$ and $\alpha$ depending only on $g$, $C_0$ and  $\| F\|_{C^2(g)}$ 
such that
$$\| \emph{tr}_{g}{g'} \|_{C^{\alpha}(g)} \le C.$$
\end{theorem}

\begin{proof} We will work locally and fix a coordinate system $(x^1, \cdots,
x^4)$ with the same properties as the one in the proof of Lemma \ref{lemmatrgg'},
with the point $p$ corresponding to $x=0$.   We will show that, with the
notation of section 3,
\begin{equation} \label{eqnholder}
| \GI{\al}{\be} \Gp{\al}{\be}(y) - \GI{\al}{\be} \Gp{\al}{\be}(x)| \le C'
R^{\alpha},
\end{equation} 
for all $x,y \in B_{R}(0)$ and $0 < R < R_0/2$ for some positive constants $\alpha$, $R_0$ and $C'$, where $B_{R}(0)$ is the ball of radius $R$ centred
at 0.  This will prove the theorem, since a short calculation shows
$$\GI{\al}{\be} \Gp{\al}{\be} = \frac{1}{2} g^{ij} g'_{ij} + O(R).$$
To prove (\ref{eqnholder}),
first note that, by a straightforward calculation,
$$(\det{G'})^2 = \frac{1}{16} \det{g'} + \eta,$$
where $\eta=\eta(x)$ is a function of the form
$$\eta \ = \sum_{a,b,c,d=1}^4 h_{abcd} G'_a G'_b G'_c G'_d,$$
where $G'_1 = \Gp{1}{1}$, $G'_2 = \Gp{2}{2}$, $G'_3 = \Gp{1}{2} + \Gp{2}{1}$
and $G'_4 = \sqrt{-1}(\Gp{1}{2} - \Gp{2}{1})$, and where the $h_{abcd}$ are
smooth functions depending only on $J$ which vanish at 0.  Note that here,
and in the sequel, we are shrinking $R_0$ whenever necessary.  Writing
$$K = \frac{1}{16} e^{2F} \det g,$$ and using the equation (\ref{eqnCY4}),
we see that
\begin{equation} \label{eqnCY5}
2 \log \det{G'} = \log (K + \eta).
\end{equation}
Define a function $\Phi$ on the space of positive definite Hermitian matrices
by $\Phi(A_{\al \overline{\be}}) = 2 \log \det (A_{\al \overline{\beta}})$.
Since $\Phi$ is concave, the
 tangent plane to the graph of $\Phi$ at a point $\Gp{\al}{\be}(y)$ lies
 above the graph of $\Phi$ and so 
$$2 \GpI{\al}{\be}(y) ( \Gp{\al}{\be}(y) - \Gp{\al}{\be}(x) ) \le \Phi(
\Gp{\al}{\be}(y) )
- \Phi( \Gp{\al}{\be}(x) ),$$
for $x,y$ in $\ov{B}_{2R}(0)$.  From (\ref{eqnCY5}),
\begin{eqnarray} \nonumber
2 \GpI{\al}{\be}(y) ( \Gp{\al}{\be}(y) - \Gp{\al}{\be}(x) ) & \le & \log
\left( 1 + \frac{K(y) - K(x) + \eta(y) - \eta(x)}{K(x) + \eta(x)} \right)
\\ \nonumber
& \le & \log( 1 + C_1 R) \\ \label{eqnconcave1}
& \le & C_1 R.
\end{eqnarray}
We now need the following elementary linear algebra lemma.

\begin{lemma} \label{lemmalinearalgebra}
Let $S(\lambda, \Lambda)$ be the set of $2 \times 2$ positive definite Hermitian
matrices with eigenvalues between $\lambda$ and $\Lambda$, with $0< \lambda
\le \Lambda$.  Then there exist a finite number of bases of unit vectors
$\{ (V_{\nu}^{(1)}, V_{\nu}^{(2)} )\}_{\nu=1}^N$ and constants $0 < \lambda^* < \Lambda^*$
depending only on $\lambda$ and $\Lambda$ such that any $A$ in $S(\lambda,\Lambda)$ can be written
$$A = \sum_{\nu=1}^N \beta_{\nu} (V_{\nu}^{(1)} \otimes \overline{V}_{\nu}^{(1)}
+ V_{\nu}^{(2)} \otimes \overline{V}_{\nu}^{(2)})$$ with $\lambda^* \le \beta_{\nu} \le
\Lambda^*$.
\end{lemma}
\begin{proof} This lemma can be proved by a straightforward modification of the argument in \cite{MoWa}.  Q.E.D.
\end{proof}
Using this lemma we see that
$$\GpI{\al}{\be}(y) = \sum_{\nu=1}^N \beta_{\nu}(y) (V_{\nu}^{(1)} \otimes \overline{V}_{\nu}^{(1)}
+ V_{\nu}^{(2)} \otimes \overline{V}_{\nu}^{(2)}),$$
for $\lambda^* < \beta_{\nu} < \Lambda^*$ where the $V^{(i)}_{\nu}$ and $\lambda^*$ and $\Lambda^*$ depend
only on the constant $C_0$.  Define
$$w_{\nu} = ( (V_{\nu}^{(1)})^{\al} \overline{(V_{\nu}^{(1)})^{\be}} + (V_{\nu}^{(2)})^{\al} \overline{(V_{\nu}^{(2)})^{\be}})\Gp{\al}{\be},$$
where $(V_{\nu}^{(i)})^{\al}$ is the $\alpha$-component of the vector $V_{\nu}^{(i)}$.
We can then rewrite (\ref{eqnconcave1}) as 
\begin{equation} \label{eqnconcave12}
\sum_{\nu=1}^N \beta_{\nu}(y) (w_{\nu}(y) - w_{\nu}(x)) \le C_1R, \qquad \textrm{for } x,y \in \ov{B}_{2R}(0).
\end{equation}
We will now use the concavity of $\Phi$ again, this time to derive a differential inequality for $w_{\nu}$.  For each $\nu$,
apply the operator $\sum_{i=1}^2 D_{\nu}^{(i)} D_{\ov{\nu}}^{(i)} = \sum_{i=1}^2 (V_{\nu}^{(i)})^{\ga} \overline{(V_{\nu}^{(i)})^{\de}}
D_{\ga} D_{\ov{\de}}$ to (\ref{eqnCY5}) to obtain:
\begin{eqnarray} \nonumber
\lefteqn{ \sum_{i=1}^2 \left( 2 \GpI{\al}{\be} D_{\nu}^{(i)} D_{\ov{\nu}}^{(i)}
\Gp{\al}{\be} - 2 \GpI{\al}{\gamma} \GpI{\sigma}{\be} (D_{\nu}^{(i)} \Gp{\sigma}{\ga})(
D_{\ov{\nu}}^{(i)} \Gp{\al}{\be}) \right) 
} \\ \label{eqnDD}
\qquad \qquad  &= & \sum_{i=1}^2 \left( \frac{D_{\nu}^{(i)} D_{\ov{\nu}}^{(i)} (K + \eta)}{K+\eta}
- \frac{| D_{\nu}^{(i)} (K+\eta)|^2}{(K+\eta)^2} \right). \qquad \qquad  \qquad \qquad \qquad 
\end{eqnarray}
Apply Lemma \ref{lemmaDG} twice to  the first term on the left hand side and
the first term on the right hand side of (\ref{eqnDD}).  Making use of the
good second term on the left hand side of (\ref{eqnDD}) we see that there
is a second order elliptic operator $L_{\nu}= a^{ij}\dd{i}\dd{j}$ with 
$C_2^{-1} |\xi|^2 \le a^{ij} \xi_i \xi_j \le C_2 |\xi|^2$ such that 
\begin{equation} \label{eqnLinequality}
L_{\nu} w_{\nu} \ge -C_3.
\end{equation}
From the inequalities (\ref{eqnconcave12}) and (\ref{eqnLinequality}) we make the following claim.

\bigskip \noindent
{\bf Claim:} There exist positive constants $\hat{C}$ and $\delta$ such that 
$$\textrm{osc}_{B_R(0)} w_{\nu} \le \hat{C} R^{\delta}, \qquad \textrm{for } 0 < R < R_0/2.$$

Of course, given this claim, we are finished, since we can then write
$$\GI{\al}{\be} \Gp{\al}{\be} = \sum_{\nu=1}^N \hat{\beta}_{\nu} w_{\nu},$$
with $\hat{\beta}_{\nu}$ smooth bounded functions depending only on $g$ and $J$ and satisfying $C_4^{-1} < \hat{\beta}_{\nu} < C_4$.  This gives (\ref{eqnholder}) and Theorem \ref{theoremholder} follows.  Q.E.D.

\bigskip \noindent
{\bf Proof of Claim} \ Although this proof can easily be extracted from \cite{Tr2} (see also \cite{Si}), we will include a sketch of the argument here for the convenience of the reader.  The key tool is the following  Harnack inequality \cite{Tr1}: if $u \ge 0$ satisfies $L u = a^{ij} \dd{i} \dd{j} u \le C_3$ with $C_2^{-1} |\xi|^2 \le a^{ij} \xi_i \xi_j \le C_2 |\xi|^2$ on $B_{2R}(0)$ then there exists $p>0$ such that
$$\left( \frac{1}{R^4} \int_{B_{R}(0)} u^p \right)^{1/p} \le C_H (\inf_{B_R(0)} u + R),$$
where the constant $C_H$ depends only on $C_2$ and $C_3$.

Set $M_{s\nu} = \sup_{B_{sR}(0)} w_{\nu}$ and $m_{s\nu} = \inf_{B_{sR}(0)} w_{\nu}$ for $s=1,2$ and apply the Harnack inequality to $(M_{2\nu} - w_{\nu})$ to obtain, for fixed $l$,
\begin{eqnarray} \nonumber
\left( \frac{1}{R^4} \int_{B_R(0)} ( \sum_{\nu \neq l} (M_{2\nu} - w_{\nu}))^p \right)^{1/p} & \le & N^{1/p} \sum_{\nu \neq l} \left( \frac{1}{R^4} \int_{B_R(0)} (M_{2 \nu} -w_{\nu})^p \right)^{1/p} \\ \nonumber
& \le & C_5 ( \sum_{\nu \neq l} (M_{2 \nu} - M_{1\nu}) + R ) \\ \nonumber \label{eqnclaim1}
& \le & C_5 ( \omega (2R) - \omega(R) + R),
\end{eqnarray}
where $\omega(sR) = \sum_{\nu=1}^N \textrm{osc}_{B_{sR}(0)} w_{\nu} = \sum_{\nu=1}^N (M_{s\nu} - m_{s\nu}).$
From (\ref{eqnconcave12}) we have
$$ \beta_l ( w_l(y) - w_l(x)) \le C_1 R + \sum_{\nu \neq l} \beta_{\nu} (w_{\nu}(x) - w_{\nu}(y)),$$
for $x,y \in \ov{B}_{2R}(0)$.  Choosing $x \rightarrow m_{2l}$ and integrating over $y \in B_R(0)$ gives
\begin{eqnarray} \nonumber
\left( \frac{1}{R^4} \int_{B_R(0)} (w_l - m_{2l})^p \right)^{1/p} & \le & C_6 R +  C_6 \left( \frac{1}{R^4} \int_{B_R(0)} ( \sum_{\nu \neq l} (M_{2 \nu} - w_{\nu} )^p \right)^{1/p} \\ \label{eqnclaim2}
& \le & C_7 (\omega(2R) - \omega(R) +R).
\end{eqnarray}
Now apply the Harnack inequality to $(M_{2l} - w_l)$ to obtain
\begin{eqnarray} \nonumber
\left( \frac{1}{R^4} \int_{B_R(0)} (M_{2l} - w_l)^p \right)^{1/p} & \le & C_8 (M_{2l} - M_{1l} + R) \\ \label{eqnclaim3}
& \le & C_8 (\omega(2R) - \omega(R) + R).
\end{eqnarray}
Combining (\ref{eqnclaim2}) and (\ref{eqnclaim3}) we see that
$$M_{2l} - m_{2l} \le C_9 (\omega(2R) - \omega(R) + R),$$
and summing in $l$ gives
$$\omega(2R) \le C_{10} (\omega(2R) - \omega(R) + R),$$
from which it follows that 
$$\omega(R) \le \left( \frac{C_{10}-1}{C_{10}} \right) \omega(2R) +  R,$$
and the claim follows by a well-known argument (see \cite{GiTr}, Chapter 8).  Q.E.D.
\end{proof}

\bigskip
\bigskip
\noindent
\setcounter{equation}{0}
{\bf 5. Higher order estimates} 
\bigskip
\addtocounter{section}{1}

In this section we will prove estimates on $g'$ and all of its derivatives given a H\"older estimate 
\begin{equation} \label{eqnestimateg'}
\| \tr{g}{g'} \|_{C^{\alpha}(g)} \le C,
\end{equation}
with $0<\al< 1$ and an estimate $g' \ge C^{-1} g$.
In light of Theorem \ref{theoremgestimate} and Theorem \ref{theoremholder} this will complete the proof of Theorem 1.

Consider the normalized almost-K\"ahler potential $\phi_0$ defined by (see section 2)
$$\Delta \phi_0 = \tr{g}{g'} - 4, \qquad \int_M \phi_0 \frac{\omega^2}{2}
=0.$$
From (\ref{eqnestimateg'}), by  the
elliptic Schauder estimates for the Laplacian we have
\begin{equation} \label{eqnphischauder}
\| \phi_0 \|_{C^{2+\alpha}(g)} \le C_0.
\end{equation}
Recall that the 1-form $a_0$ satisfies
$$\omega' = \omega - \frac{1}{2} d (J d\phi_0) + da_0.$$
Without loss of generality, we can assume that $a_0$ is $L^2$ orthogonal
to the harmonic 1-forms.  Then since $a_0$ satisfies the uniformly elliptic system (\ref{eqnellipticsystem1.5}) for $s=0$ and is orthogonal to its kernel we can use
 (\ref{eqnphischauder}) and the
Schauder  elliptic estimates  to obtain
$$\| a_0 \|_{C^{2+ \alpha}(g)} \le C_1.$$
It follows that $\| g' \|_{C^{\alpha}(g)} \le C_2$.
Differentiating the Calabi-Yau
equation (\ref{eqnCY4}), we see that
\begin{equation} \label{eqnelliptic}
g'^{ij} \dd{i} \dd{j} (\dd{k} \phi_0) + \{ \textrm{lower
order terms} \} = g^{ij} \dd{k} g_{ij} + 2\dd{k}F,
\end{equation}
where the lower order terms may contain up to two derivatives of $\phi_0$ or
$a_0$.  Since the coefficients of this elliptic equation are in $C^{\alpha}$
we can apply the standard Schauder estimates again to obtain $$\| \phi_0 \|_{C^{3+
\alpha}(g)} \le C_3.$$
From (\ref{eqnellipticsystem1.5}) we then obtain
$$\| a_0 \|_{C^{3+ \alpha}(g)} \le C_4.$$  The rest of the higher order
estimates follow from (\ref{eqnelliptic}),  (\ref{eqnellipticsystem1.5}) and
a bootstrapping argument.  This completes the proof of Theorem 1.  Q.E.D.

\setcounter{equation}{0}
\setcounter{lemma}{0}
\addtocounter{section}{1}
\bigskip
\bigskip
\noindent
{\bf 6. Proof of Theorem 2: the case $\mathbf{b^+(M)=1}$}
\bigskip

For this section we assume $b^+(M)=1$.  Consider the equation
\begin{equation} \label{eqncontinuity}
{\omega'_t}^2  =  e^{tF + c_t} \omega^2,
\end{equation}
where $c_t$ is the constant given by $c_t = \log ( \int_M \omega^2 / \int_M e^{tF} \omega^2)$, and where $\omega'_t$ is required to be cohomologous to $\omega$ and compatible with $J$.
Let 
$$T = \{ t' \in [0,1] \ | \ \exists \textrm{ smooth } \omega'_t \textrm{ solving } (\ref{eqncontinuity}) \textrm{ for } t \in [0,t'] \}.$$
Clearly $0 \in T$.  We will show that $T$ is both open and closed in $[0,1]$.  This will prove Theorem 2.  Note that if $\omega'_t$ is in $C^{\alpha}$ then by the estimates of section 5 it is smooth.  

We show now that $T$ is open.  Fix $t_0$ in $T$.  We will show that (\ref{eqncontinuity}) can be solved for $t$ in an open neighbourhood containing $t_0$.  Fix $\ot = \omega'_{t_0}$.  Then solving (\ref{eqncontinuity}) near $t_0$  is equivalent to solving
$$\log \frac{{\omega'_t}^2}{\ot^2} - (t-t_0)F - (c_t - c_{t_0}) = 0,$$
for $t$ close to $t_0$.  

Let $\Lambda^{k,s+\alpha}$ be the space of $k$-forms in $C^{s+\alpha}$, and let $W^{\alpha} \subset \Lambda_{\ot}^+$ be the space of self-dual two forms $\gamma$ in  $C^{\alpha}$ satisfying $\int_M \exp\left(\frac{2 \gamma \wedge \tilde{\omega}}{\tilde{\omega}^2}\right) \tilde{\omega}^2 = \int_M \tilde{\omega}^2$.  Then define a map
$$\Phi: \Lambda^{1,1+\alpha} \times \mathbb{R} \rightarrow W^{\alpha},$$
by
$$\Phi(b, t) = \left( \log \frac{(\ot+ db)^2}{\ot^2} - (t-t_0)F - \hat{c} \right) \frac{\ot}{2} +  \PP{}{}{}{} db,$$
where $$\hat{c}(b,t) =  \log \left( \int_M e^{-(t-t_0)F} (\ot + db)^2 \right)
- \log \int_M \omega^2.$$   Note that if we can find $b=b(t)$ solving $\Phi(b,t)=0$ for $t$ close to $t_0$, then this would imply $\hat{c} = c_{t} - c_{t_0}$ and complete the openness argument.  

Since $b^+(M)=1$, the space $\mathcal{H}^+_{\ot}$ of harmonic
self-dual 2-forms with respect to $\ot$ is spanned by $\ot$.  Notice that the tangent space to $W^{\alpha}$ at $\Phi(0,t_0)$ is equal to $(\mathcal{H}_{\ot}^+)^{\perp} \cap \Lambda^{2,\alpha}$, where $(\mathcal{H}_{\ot}^+)^{\perp} \subset \Lambda_{\ot}^+$ is the space of self-dual 2-forms which are $L^2(\ot)$ orthogonal to $\mathcal{H}_{\ot}^+$.  Then the derivative of $\Phi$ in the $b$-variable at $(0,t_0)$ is a map
$$(D_1\Phi)_{(0,t_0)}: \Lambda^{1,1+\alpha} \rightarrow (\mathcal{H}_{\ot}^+)^{\perp} \cap \Lambda^{2,\alpha}$$
given by
$$(D_1 \Phi)_{(0,t_0)} (\beta) = d_{\ot}^+ \beta.$$
It is well known (see \cite{DoKr},
for example) that this map is surjective and hence openness follows by the implicit function theorem.

\bigskip

We now need to prove closedness under the assumption that the Nijenhuis tensor is small in the $L^1$ sense.  Note that from the discussion in section 2, since $|\nabla N(J)|$ is uniformly bounded in terms of the curvature of $g$, if the Nijenhuis tensor is small in the $L^1$ norm, it is small in the $C^0$ norm, and hence also in the $L^p$ norm for any $p>1$.   It will be convenient (and sufficient) for us to prove Theorem 2 under the assumption that $N(J)$ is small in some $L^p$ norm, where $p$ will be a fixed constant strictly larger than 2.

We will use the following lemma.

\begin{lemma} \label{lemmaphi1}
Let $\omega' = \omega - \frac{1}{2} d(J d\phi_1) + da_1$ be a solution of the Calabi-Yau equation (\ref{eqnCY2}).  Suppose that for some constants $p>2$ and $B$,
\begin{equation} \label{eqnlemmaphi1}
\left( \int_M \left| \frac{da_1 \wedge \omega}{\omega^2} \right|^p \omega^2 \right)^{1/p} \le B.
\end{equation}
Then there exists a constant $C'$ depending only on $g$, $p$, $B$ and $\sup_M|F|$ such that
$$\sup_M \phi_1 - \inf_M \phi_1 \le C'.$$
\end{lemma}
\begin{proof} This is a modification of Yau's  well-known Moser iteration argument.
For ease of notation,  write $\phi=\phi_1$.  Assume that $\int_M \phi \, \omega^2=0$. For $\l \ge 0,$
\begin{eqnarray*}
C_0 \int_M | \phi|^{l+1} \omega^2 & \ge & \int_M \phi | \phi|^l (\omega^2 - \omega'^2) \\
& = & \frac{1}{2} \int_M \phi | \phi |^l d (J d\phi) \wedge (\omega + \omega') - \int_M \phi | \phi|^l da_1 \wedge \omega \\
& = & - \frac{(l+1)}{2} \int_M | \phi |^l d\phi \wedge Jd \phi \wedge (\omega + \omega') - \int_M \phi | \phi|^l da_1 \wedge \omega \\
& = & - \frac{(l+1)}{2(l/2 + 1)^2} \int_M d (\phi | \phi |^{l/2}) \wedge Jd (\phi | \phi |^{l/2}) \wedge (\omega + \omega') \\
&& \mbox{} - \int_M \phi | \phi|^l da_1 \wedge \omega \\
& \ge & \frac{(l+1)}{4(l/2+1)^2} \int_M | \nabla ( \phi | \phi|^{l/2})|^2 \omega^2 \\
&& \mbox{} -  \left( \int_M | \phi |^{q(l+1)} \omega^2  \right)^{1/q} \left( \int_M \left| \frac{da_1 \wedge \omega}{\omega^2} \right|^p \omega^2 \right)^{1/p},
\end{eqnarray*}
for $q$ satisfying $1/p + 1/q =1$.  Setting $l=0$ we see that since $q<2$,
\begin{eqnarray*}
\int_M | \nabla \phi |^2 & \le & C_1 \left( \left( \int_M | \phi |^q  \right)^{1/q} + \int_M | \phi | \, \right) 
 \le  C_2\left( \left(\int_M | \phi |^2 
 \right)^{1/2} + 1\right),
\end{eqnarray*}
where we have omitted the volume form $\omega^2$.
Since $\int_M \phi  = 0$ we obtain $\| \phi \|_{L^2} \le C_3$ from the Poincar\'e inequality.  

We have
for general $l$,
\begin{equation} \label{eqnreversesobolev}
\int_M | \nabla ( \phi | \phi|^{l/2})|^2  \le C_4(l+2) \max \left\{ 1, \int_M | \phi |^{l+2} , \left(  \int_M | \phi |^{q(l+1)}  \right)^{1/q} \right\}.
\end{equation}
The Sobolev inequality  gives
$$ \left( \int_M |u|^4  \right)^{1/2} \le C_5 \left( \int_M | \nabla u|^2   + \int_M u^2  \right),$$
for functions $u$ on $M$. 
Set $r=l+2 \ge 2$.   Applying the Sobolev inequality  
 to $u= \phi | \phi|^{l/2}$, making use of (\ref{eqnreversesobolev}) and raising to the power $1/r$ gives 
$$\| \phi \|_{L^{2r}} \le C_6^{1/r} r^{1/r} \max \{ 1, \| \phi \|_{L^r}, 
\| \phi \|^{(r-1)/r}_{L^{q(r-1)}} \}.$$
Setting $r=2$ we obtain $\| \phi \|_{L^4} \le C_7$.  For general $r$ we use that fact that $\| \phi \|_{L^a} \le C_8 \| \phi \|_{L^b}$ whenever $a \le b$ to see that
$$\| \phi \|_{L^{2r}} \le C_9^{1/r} r^{1/r} \max \{ 1, 
\| \phi \|_{L^{qr}} \}.$$
By successively replacing $r$ by $\sigma r$ for $\sigma = 2/q>1$ we see that for all $k=0, 1,2, \ldots$,
$$\| \phi \|_{L^{2r\sigma^k}} \le C_9^{( \frac{1}{r} \sum_{i=0}^k \sigma^{-i})} r^{(\frac{1}{r} \sum_{i=0}^k \sigma^{-i})} \sigma^{(\frac{1}{r} \sum_{i=1}^k i\sigma^{-i})} \max \{ 1 , \| \phi \|_{L^{qr}} \}.$$
Set $r = 2 $ and let $k \rightarrow \infty$.  This gives a bound $\| \phi \|_{C^0} \le C_{10} \max \{ 1, \| \phi \|_{L^{2q}} \}$, which is uniformly bounded since $q<2$. This completes the proof of the lemma.  Q.E.D.
\end{proof}

\noindent
{\bf Remark 6.1} \  In a private discussion, Donaldson made the following
surprising observation:  the almost-K\"ahler potential $\phi_{1/2}$ is  uniformly bounded, without any assumption on $N(J)$. This can be proved using  a Moser iteration argument and equation (\ref{eqnOmega}).

\bigskip

It is now not difficult to complete the proof of Theorem 2.  We suppose that we have a solution of (\ref{eqncontinuity}) for $t \in [0,t_0)$ for some $t_0 \in [0,1)$. 
  We require uniform estimates on $\omega_t'$ and all its derivatives and by Theorem 1, it is sufficient to obtain a uniform estimate of $\phi_1$.  We have the following claim. 

\bigskip
\noindent
{\bf Claim:} \ Let $p>2$.  There exists $\epsilon>0$ depending only on $p$,
$g$ and $\| F \|_{C^2(g)}$ such that if $\| N(J) \|_{L^{p}} < \epsilon$ then for $t \in [0,t_0)$, $a_1 = a_1(t)$ satisfies 
\begin{equation} \label{eqnclaim}
\left( \int_M \left| \frac{da_1 \wedge \omega}{\omega^2} \right|^p \omega^2 \right)^{1/p} <1.
\end{equation}
 
\bigskip
\noindent
{\bf Proof of Claim} \ At $t=0$ we  have $da_1=0$.  Suppose that the claim is false.  Then there is a first time $t' \in [0,t_0)$ with 
$$\left( \int_M \left| \frac{da_1 \wedge \omega}{\omega^2} \right|^p \omega^2 \right)^{1/p}  = 1.$$  
Then it follows from Lemma \ref{lemmaphi1}  that we have a H\"older estimate on $\omega'$  at $t=t'$.  Now the $L^p$ $\emph{a priori}$ estimates for  the elliptic system (\ref{eqnellipticsystem1.5}) with $s=1$ give
$$ \left( \int_M \left| \frac{da_1 \wedge \omega}{\omega^2} \right|^p \omega^2 \right)^{1/p} \le K \| N(J) \|_{L^{p}},$$
for some uniform constant $K$.  Picking $\epsilon = 1/2K$ gives a contradiction and proves the claim.  Q.E.D.

\bigskip

Then the first part of Theorem 2 follows from this claim and the previous
lemma.

\bigskip

\setcounter{equation}{0}
\setcounter{lemma}{0}
\addtocounter{section}{1}
\pagebreak[3]
\noindent
{\bf 7. Proof of Theorem 2: the case  $\mathbf{b^+(M)>1}$}

\bigskip

Suppose that $b^+(M)=r+1$.
We begin with the openness argument.  For convenience,  assume that 
 $\omega$ has been scaled so that $\int_M
\omega^2=1$.  We wish to solve the equation
\begin{equation} \label{eqncontinuity2}
{\omega'_t}^2 = e^{tF+c_t} \omega^2,
\end{equation}
with $c_t = - \log ( \int_M e^{tF} \omega^2)$, for $\omega'_t$ satisfying  $\int_M \omega'_t \wedge \omega >0$ and 
$[\omega'_t] \in H^+_{\omega}$.
As in section 6, we suppose that there is a solution $\ot = \omega'_{t_0}$ at $t=t_0$ and show that  (\ref{eqncontinuity2}) can be solved for $t$ close to $t_0$.

Let $\chi_1, \ldots, \chi_r$ and $\tilde{\chi}_1, \ldots, \tilde{\chi}_r$
 be self-dual harmonic 2-forms with respect to $\omega$ and $\ot$ respectively such that $\{ \omega, \chi_1, \ldots, \chi_r \}$ and $\{ \ot, \tilde{\chi}_1,
\ldots, \tilde{\chi}_r \}$ are $L^2$ orthonormal bases for $\mathcal{H}^+_{\omega}$
and $\mathcal{H}^+_{\ot}$.
Let $\Lambda^{1, 1+\alpha}$ and $W^{\alpha}$ be as in section 6.
Consider the operator $\Phi : \Lambda^{1, 1+\alpha} \times \mathbb{R}^r \times \mathbb{R} \rightarrow W^{\alpha}$ given by
$$\Phi(b,\un{s},  t) = \left( \log \frac{(\ot  + \sum_{i=1}^r s_i \chi_i+ db)^2}{\ot^2} - (t-t_0)F - \hat{c} \right) \frac{\ot}{2} +  \PP{}{}{}{} ( \sum_{i=1}^r s_i \chi_i + db),$$
where $$\hat{c}(b,\un{s}, t) =  \log \left( \int_M e^{-(t-t_0)F} (\ot + \sum_{i=1}^r s_i \chi_i+ db)^2 \right).$$
We have a solution $\Phi(0,\un{0}, t_0)=0$ and 
if we can find $b$ and $\un{s}$ depending on $t$ solving $\Phi(b, \un{s}, t) =0$ for $t$ near $t_0$, then after rescaling we would have our desired
solution.
Write $\Pi_{\langle \tilde{\chi}_1, \ldots, \tilde{\chi}_r \rangle}$ for the $L^2(\ot)$ projection  onto the space spanned by $\tilde{\chi}_1, \ldots, \tilde{\chi}_r$.  Define  $\Psi_1 :  \Lambda^{1, 1+\alpha} \times \mathbb{R}^r \times \mathbb{R} \rightarrow W^{\alpha} \cap (\langle \tilde{\chi}_1, \ldots, \tilde{\chi}_r \rangle)^{\perp}$ by
$$\Psi_1 (b, \un{s}, t)  =  (1 - \Pi_{\langle \tilde{\chi}_1, \ldots, \tilde{\chi}_r \rangle}) \Phi(b, \un{s}, t).$$
The derivative $(D_1 \Psi_1)_{(0,\un{0}, t_0)} : \Lambda^{1,1+\al} \rightarrow (\mathcal{H}_{\ot}^+)^{\perp} \cap \Lambda^{2, \alpha} $ is given by
$$(D_1 \Psi_1)_{(0,\un{0}, t_0)} (\beta) = d_{\ot}^+\beta.$$
This map is surjective and so by the implicit function theorem, given $(\un{s}, t)$ near $(\un{0},t_0) \in \mathbb{R}^r \times \mathbb{R}$ there exists $b= b(\un{s}, t)$ solving $\Psi_1 = 0$.  Now define a map $\Psi_2 : \mathbb{R}^r \times \mathbb{R} \rightarrow \mathbb{R}^r$ in a neighbourhood of $(\un{0}, t_0)$ by
$$\Psi_2 (\un{s}, t) = \Pi_{\langle \tilde{\chi}_1, \ldots, \tilde{\chi}_r \rangle} \Phi(b(\un{s}, t), \un{s}, t),$$
where we are identifying $\mathbb{R}^r$ and the space spanned by the $\tilde{\chi}_i$.  Calculate
$$\left( (D_1 \Psi_2 \right)_{(\un{0},t_0)} )_{ij} = \int_M \left \langle \frac{1}{2} (1+*_{\ot})
\left(\chi_j +  d \left( \frac{\partial b}{\partial s_j}(\un{0},t_0) \right) \right), \tilde{\chi}_i \right \rangle_{\ot} \frac{\ot^2}{2}
= \int_M \chi_j \wedge \tilde{\chi}_i,$$
which is invertible.
Applying the implicit function theorem again we find $\un{s} = \un{s}(t)$ solving $\Psi_2(\un{s}(t), t)=0$ and hence
$$\Phi(b(\un{s}(t), t), \un{s}(t), t) =0,$$
for $t$ close to $t_0$.  This completes the proof of openness.  

We now turn to the question of closedness. As discussed in section 6, we may assume that that $N(J)$ is small in the $L^p$ sense for any fixed $p$.  Assume that we have a solution of (\ref{eqncontinuity2})
with $\int_M \omega'_t \wedge
\omega>0$ and $[\omega'_t] \in H^+_{\omega}$ on some maximal  interval
$[0,t_0)$.  Write $\omega'= \omega'_t$  and
define $s_i$ by $[\omega'] = [\omega] + \sum_{i=0}^r s_i [\chi_i]$ for $\chi_0
=\omega$ and $\chi_1, \ldots, \chi_r$
as above.  Notice that by squaring both sides of this equation we see that
the $s_i$ are bounded.  Define $\phi_0$ and $\phi_1$ by 
\begin{eqnarray} \label{eqnnewphi0}
\frac{1}{4} \Delta \phi_0 & = & \frac{\omega \wedge \omega'}{\omega^2} - \int_M \omega' \wedge \omega \\ \label{eqnnewphi1}
\frac{1}{4} \Delta' \phi_1 &= & \int_M \omega' \wedge \omega - \frac{\omega \wedge \omega'}{\omega'^2},
\end{eqnarray}
where we recall that $\int_M \omega^2 = \int_M \omega'^2=1$.
Let us first assume that $\phi_1$ is uniformly bounded.  Then
since $\int_M \omega' \wedge \omega = 1 + s_0>0$ is uniformly bounded from
above, Theorem \ref{theoremgestimate} still holds with essentially the same proof.  Notice that the  bound on $\tr{g}{g'}$ implies a uniform positive
lower bound
for $\int_M \omega' \wedge \omega$.  No changes are necessary for section 4.  For the higher order estimates, we argue as follows.   Define $a_0$ by $d^*a_0=0$ and
$$\omega' = \omega + \sum_{i=0}^r s_i \chi_i - \frac{1}{2} d(Jd\phi_0) + da_0.$$
Then $a_0$ satisfies the equations
\begin{eqnarray*}
 da_0 \wedge \omega & = & - \sum_{i=1}^r s_i \chi_i \wedge \omega \\
\PP{}{}{}{}{} da_0 & =& - \sum_{i=1}^r s_i \PP{}{}{}{} \chi_i + \frac{1}{4} (\dd{i} \J{j}{k} - \dd{j} \J{i}{k}) (\dd{k} \phi_0) \, dx^i \wedge dx^j,
\end{eqnarray*}
and the arguments of section 5 follow in just the same way as before.

We will now show that $\phi_1$ can be bounded if $N(J)$ is small in the $L^p$ norm for $p>2$.  
Define $a_1$ by 
$$\omega' = \omega + \sum_{i=0}^r s_i \chi_i - \frac{1}{2} d(J d\phi_1) + da_1,$$
and $d^*_1a_1=0$ where we are using the subscript 1 to denote the metric $\omega'$.
For ease of notation, set $$\zeta = \frac{1}{4} (\dd{i} \J{j}{k} - \dd{j} \J{i}{k}) (\dd{k} \phi_1) \, dx^i \wedge dx^j.$$ The 1-form $a_1$ satisfies
\begin{eqnarray*}
da_1 \wedge \omega' & = & - s_0 \omega'^2 -  \sum s_i \chi_i \wedge \omega' \\
\PP{}{}{}{} d a_1 & = &   \zeta  - \PP{}{}{}{}(\sum s_i \chi_i),
\end{eqnarray*}
where, here and from now on, we are always summing $i$ from $0$ to $r$.  This equation
can  be rewritten as
\begin{equation} \nonumber \label{eqndplusa1}
d^+_1 a_1 = \zeta - s_0 \omega' - \frac{1}{2} (1+*_1) \sum s_i \chi_i.
\end{equation}
Write $\Pi$ for the $L^2(\omega')$ projection onto the space $\mathcal{H}_{\omega'}^+$ of self-dual harmonic forms with respect to $\omega'$.  Then we see that
\begin{eqnarray} \label{eqnzeta1}
(1-\Pi) \zeta & = & d_1^+ a_1 + (1- \Pi) \frac{1}{2} (1+*_1)   \sum s_i \chi_i \\ \label{eqnzeta2}
 \Pi \zeta & =  &  s_0\omega' + \Pi \sum s_i \chi_i = \Pi \sum \tilde{s}_i
 \chi_i,
\end{eqnarray}
for $\tilde{s}_0 = s_0(2+s_0)$ and $\tilde{s}_i = s_i (1+s_0)$ for $i\ge 1$.
Now Lemma \ref{lemmaphi1} holds as before if we replace (\ref{eqnlemmaphi1}) by the inequality
$$ \left( \int_M \left| \frac{da_1 \wedge \omega}{\omega^2} \right|^p \omega^2 \right)^{\frac{1}{p}} + \left( \int_M \left| \frac{ \sum s_i \chi_i \wedge \omega}{\omega^2} \right|^p \omega^2 \right)^{\frac{1}{p}} \le B,$$
where we are making use of the fact that $(da_1 + \sum s_i \chi_i )\wedge \omega' = - s_0 \omega'^2$ is bounded.  We can now replace the inequality (\ref{eqnclaim}) in the Claim by
\begin{equation} \label{eqnnewclaim}
 \left( \int_M \left| \frac{da_1 \wedge \omega}{\omega^2} \right|^p \omega^2 \right)^{\frac{1}{p}} + \left( \int_M \left| \frac{  \sum s_i \chi_i \wedge \omega}{\omega^2} \right|^p \omega^2 \right)^{\frac{1}{p}} < 1.
\end{equation}
Indeed, arguing for a contradiction as in the proof of the claim, we suppose that we have equality in (\ref{eqnnewclaim}).  Then $\omega'$ and $\omega$ are uniformly equivalent and we can essentially ignore the fact that they define different norms.  Writing $C$ for a uniform constant which may change from inequality to inequality we have $\| \zeta \|_{L^p} \le C \| N(J) \|_{L^p}$ from which it follows that $ \| \Pi \zeta \|_{L^2} \le C \| N(J) \|_{L^p}$. 
Then we see from
 (\ref{eqnzeta2}) that $|\tilde{s}_i| \le C \| N(J) \|_{L^p}$. Hence $|s_i|
 \le C \| N(J) \|_{L^p}$ and 
\begin{equation}  \label{eqnestimatesum}
\left\| \sum s_i \chi_i \right\|_{L^p} \le C \| N(J) \|_{L^p}. 
\end{equation}
But  we also have  $\| \Pi
\zeta \|_{L^p} \le C \| N(J) \|_{L^p}$ and hence  $\| (1-\Pi) \zeta \|_{L^p} \le C \| N(J) \|_{L^p}$.  Then from (\ref{eqnzeta1})
 and  the elliptic $L^p$ estimates we have 
\begin{equation} \label{eqnestimatea1}
\| da_1 \|_{L^p} \le C \| N(J) \|_{L^p}.
\end{equation}
Choosing $\epsilon$ sufficiently small, we obtain the contradiction from (\ref{eqnestimatesum}) and (\ref{eqnestimatea1}).  Q.E.D.

\bigskip
\noindent
{\bf Acknowledgements.} \ The author is very grateful to:
Simon Donaldson  for suggesting this problem and for many subsequent helpful, insightful and encouraging discussions; S.-T. Yau whose original paper \cite{Ya} made this work possible and whose lectures and discussions at Harvard University have been an invaluable source of inspiration;
the author's former advisor D.H. Phong for his continued support and encouragement;  Richard Thomas for some very useful conversations and for his help in  arranging the author's year-long visit to Imperial College;
 Xiuxiong Chen, Joel Fine, Mark Haskins, Tom Mrowka, Jian Song and Valentino Tosatti for some helpful discussions.


\begin{thebibliography}{99}
\bibitem[Au]{Au} Aubin, T. {\em \'{E}quations du type {M}onge-{A}mp\`ere sur les vari\'et\'es k\"ahl\'eriennes compactes}, Bull. Sci. Math. (2) {\bf
102} (1978), no. 1, 63--95
\bibitem[BaKo1]{BaKo1} Bando, S. and Kobayashi, R. {\em Ricci-flat K\"ahler
metrics on affine algebraic manifolds}, In {\em Geometry and Analysis on
Manifolds}, Volume 1339 of {\em Lecture Notes in Mathematics}, 20--31, Springer-Verlag,
1988
\bibitem[BaKo2]{BaKo2} Bando, S. and Kobayashi, R. {\em Ricci-flat K\"ahler
metrics on affine algebraic manifolds. II}, Mathematische Annalen {\bf 287} (1990),
175--180 
\bibitem[Ca]{Ca} Calabi, E. {\em The space of K\"ahler metrics}, In {\em
Proceedings of the International Congress of Mathematicians, Amsterdam, 1954},
Vol. 2, 206--207, North-Holland, Amsterdam, 1956
\bibitem[De]{De} Delano\"e, P. {\em Sur l'analogue presque-complexe de l'\'equation de Calabi-Yau}, Osaka J. Math. {\bf 33} (1996), no. 4, 829--846
\bibitem[Do]{Do} Donaldson, S.K. {\em Two-forms on four-manifolds and elliptic equations}, in preparation
\bibitem[DoKr]{DoKr} Donaldson, S.K. and Kronheimer, P.B. {\em The geometry
of four-manifolds}, Clarendon Press, Oxford, 1990
\bibitem[Ev]{Ev} Evans, L.C. {\em Classical solutions of fully nonlinear, convex,
second order elliptic equations}, Comm. Pure Appl. Math {\bf 25} (1982),
333--363
\bibitem[FeGoGr]{FeGoGr} Fern\'andez, M., Gotay, M. and Gray, A. {\em Compact
parallelizable four-dimensional symplectic and complex manifolds}, Proc.
Amer. Math. Soc. {\bf 103} (1988), no. 4, 1209--1212
\bibitem[GiTr]{GiTr} Gilbarg, D. and Trudinger, N.S. {\em Elliptic partial
differential equations of second order}, Springer-Verlag, Berlin, 1977
\bibitem[Gr]{Gr} Gromov, M. {\em Pseudo holomorphic curves in symplectic manifolds}, Invent. Math. {\bf 82} (1985), 307--347
\bibitem[Ko]{Ko} Kolodziej, S. {\em
The complex Monge-Amp\`ere equation}, Acta Math. {\bf 180} (1998), no. 1, 69--117
\bibitem[Kr]{Kr} Krylov, N.V. {\em Boundedly nonhomogeneous elliptic and parabolic
equations}, Izvestia Akad. Nauk. SSSR {\bf 46} (1982), 487--523.  English
translation in Math. USSR Izv. {\bf 20} (1983), no. 3, 459--492
\bibitem[Jo]{Jo} Joyce, D. {\em Compact manifolds with special holonomy},
Oxford Mathematical Monographs, Oxford University Press, Oxford, 2000
\bibitem[MoWa]{MoWa} Motzkin, I. and Wasow, W. {\em On the approximation
of linear elliptic differential equations by difference equations with positive
coefficients}, J. Math. Phys. {\bf 31} (1952), 253--259
\bibitem[Si]{Si} Siu, Y.-T. {\em Lectures on Hermitian-Einstein metrics for
stable bundles and K\"ahler-Einstein metrics}, DMV Seminar, Volume 8, Birkhauser
Verlag, Basel, 1987
\bibitem[Th]{Th} Thurston, W. {\em Some simple examples of symplectic manifolds},
Proc. Amer. Math. Soc. {\bf 55} (1976), 467--468
\bibitem[TiYa1]{TiYa1} Tian, G. and Yau, S.-T. {\em Complete K\"ahler manifolds
with zero Ricci curvature I}, J. Amer. Math. Soc. {\bf 3} (1990), no. 3, 579--609
\bibitem[TiYa2]{TiYa2} Tian, G. and Yau, S.-T. {\em Complete K\"ahler manifolds
with zero Ricci curvature II}, Invent. Math. {\bf 106} (1991), no. 1, 27--60
\bibitem[Tr1]{Tr1} Trudinger, N.S. {\em Local estimates for subsolutions
and supersolutions of general second order elliptic quasilinear equations},
Invent. Math. {\bf 61} (1980), 67--79
\bibitem[Tr2]{Tr2} Trudinger, N.S. {\em Fully nonlinear, uniformly elliptic
 equations under natural structure conditions}, Trans. Amer. Math. Soc. {\bf
278} (1983), no. 2, 751--769
\bibitem[Ya]{Ya} Yau, S.-T. {\em On the Ricci curvature of a compact K\"ahler
manifold and the complex Monge-Amp\`ere equation, I}, Comm. Pure Appl. Math.
{\bf 31} (1978), no.3, 339--411
\end{thebibliography}
\end{document}